\documentclass[12pt]{article}
\usepackage{}
\usepackage{amsfonts}
\usepackage{latexsym}
\usepackage{color, amsmath,amssymb, amsfonts, amstext,amsthm, latexsym, dsfont}
\usepackage{epic,eepic}
\usepackage{graphicx}
\usepackage{longtable}
\usepackage{float}
\usepackage[colorlinks=false, linktocpage=true]{hyperref}
\usepackage{algorithm}
\usepackage{algorithmic}
\usepackage{float}
\usepackage{lipsum}
\makeatletter
\newenvironment{breakablealgorithm}
{
	\begin{center}
		\refstepcounter{algorithm}
		\hrule height.8pt depth0pt \kern2pt
		\renewcommand{\caption}[2][\relax]{
			{\raggedright\textbf{\ALG@name~\thealgorithm} ##2\par}%
			\ifx\relax##1\relax 
			\addcontentsline{loa}{algorithm}{\protect\numberline{\thealgorithm}##2}%
			\else 
			\addcontentsline{loa}{algorithm}{\protect\numberline{\thealgorithm}##1}%
			\fi
			\kern2pt\hrule\kern2pt
		}
	}{
		\kern2pt\hrule\relax
	\end{center}
}
\makeatother

\numberwithin{equation}{section} \topmargin -0.4in

\renewcommand{\vec}[1]{\boldsymbol{#1}} 
\title{Identifying stochastic governing equations from data of the most probable transition trajectories
\footnotetext{$^*$ Corresponding author:  duan@iit.edu}}
\author{Jian Ren$^a$ and Jinqiao Duan$^{b*}$   \\
  \small $^a$ College of Mathematics and Information Science, Zhengzhou University of Light Industry, \\
 \small Zhengzhou, 450002,  China  \\
 \small $^b$ Department of Applied Mathematics, College of Computing, Illinois Institute of Technology, \\
 \small Chicago,  IL 60616,  USA }

\begin{document}
\date{\today }

\maketitle

\pagestyle{plain}
\begin{abstract}
Extracting   governing stochastic differential equation models from elusive data is crucial to understand and forecast dynamics for complex systems. We devise a method to extract the drift term and estimate the diffusion coefficient of a governing stochastic dynamical system, from its time-series data of the most probable transition trajectory. By the Onsager-Machlup theory, the most probable transition trajectory satisfies the corresponding Euler-Lagrange equation, which is a second order deterministic ordinary differential equation involving the drift term and diffusion coefficient. We first estimate the coefficients of the Euler-Lagrange equation based on the data of the most probable trajectory, and then we calculate the drift and diffusion coefficients of the governing stochastic dynamical system. These two steps involve sparse regression and optimization.  Finally, we illustrate our method with an example and some discussions.




\end{abstract}
\section{Introduction}

Discovering relationships between the observable features and responses is a   regression  problem of supervised learning. We often formulate the response variable as a deterministic expression of the features, e.g., in linear, polynomial, trigonometric function basis , together with a random error with zero mean and  independent of the features  \cite{Hastie}.
Complex dynamical systems are usually governed by a ordinary or partial differential equation. With  time-series observations, Brunton et al. and Rudy et al. \cite{Brunton, Rudy} devised a method to find  the governing equations, by a sparse regression algorithm.

To determine a governing equation,  we face two significant issues --- selecting appropriate measurable features as basis and sparsing the coefficients of the features \cite{Yair}. There is no universal basis that performs well for all the cases. However, based on the theories of Taylor expansion and Fourier series   expansion, polynomials and trigonometric functions   are naturally  selected as features \cite{Brunton, Rudy}.
The remaining sparse regression is aimed at making some or more  coefficients zero to eliminate the uncorrelated features and prevent overfitting caused by minimizing least-square loss. The $L_0$ norm, which counts the number of the non-zero coefficients, is naturally added to the loss function as a penalty. However, it makes the loss function nonconvex and the optimization problem becomes intractable to solve \cite{Natraajan}. Consequently, the penalty is then relaxed to the less shrinkage but resolvable $L_2$ (Ridge) regression or the $L_1$ norm (Lasso) \cite{ChenSS, Tibshirani}. Comparatively, $L_1$ is sparser than $L_2$. In order to be closer proximate the sparsity of the $L_0$ norm, some properties and iterative algorithms of the nonconvex $L_{1/2}$ norm are studied \cite{Chen, Xu2010, Zeng}.

The observations are unavoidably contaminated by external noise, no matter the data are for the response variables \cite{Xu2010} or for time-series trajectories of processes \cite{Brunton, Rudy}. Facing the undesirable noise,  one way \cite{Brunton, Rudy} is to filter it, and employ a  total variation regularization to denoise the derivatives and then construct a data matrix. However, stronger noise would generally result in greater deviation for the sparse solution,   a better penalty is employed or a logical sparse regression framework is proposed, for example, Sparse Identification of Nonlinear Dynamics (SINDy) \cite{Brunton, Rudy}. It  incessantly  vanish the estimators with absolute values smaller than a threshold, and reestimate the non-zero coefficients until a stable solution is  reached. The threshold is selected by cross validation as a value minimizing the loss  function.

In this present paper, we present a method  to extract a governing stochastic differential equation,  with the drift term  and a diffusion coefficient, from a time-series observation of  its most probable transition trajectory.
Here, the   diffusion term describes      external random influences.
For these stochastic systems, the   transition pathways or the most probable transition pathways  may be observed  in  some systems. For example, the  metabolic pathways on a prokaryotic genome  \cite{Iwasaki},  and the primary absorbance and fluorescence signals for gene expression \cite{Stefan}.
Technically, the one-dimensional diffusion has also been predicted by extension (time) domain averages, which match the dominant path shape for transition pathways \cite{Hoffer}.
The most probable transition trajectory, as substitution of all the pathways in the state space, has been used to approximate the transition probability density function \cite{Psaros}.
The most probable  pathways between different brain regions have been simulated by a probabilistic fiber tractography algorithm, called ConTrack algorithm, from measurable pathways \cite{ConTrack}.
Theoretically,  certain statistical characterizations (e.g., mean exit time and escape probability)  for  a stochastic differential equation satisfy  determined differential equations \cite{WuFuDuan}.
The most probable transition trajectory in a finite time interval with fixed boundary points, of a diffusion process  (e.g., a solution of a stochastic differential equation) can also be characterized by a second-order ordinary differential equation (ODE) with boundary conditions \cite{Durr}.
This is attained by expressing the transition probability of the diffusion process as an integral, whose integrand is called Onsager-Machlup (OM) function.
And then a variational principle reveals that the OM function satisfies the Euler-Lagrange equation, which furthermore leads to a desired second-order ODE. The corresponding boundary value problem is   numerically solved by a shooting method \cite{Kress, Kutz}.
The OM function and most probable transition pathway for  a jump-diffusion process  has been studied recently by Chao and Duan  \cite{ChaoDuan}.

This paper is arranged as follows.  We first present the background of OM function for a diffusion system, and the Euler-Lagrange equation for the most probable transition trajectory, i.e.,  the   second-order ODE with boundary conditions,  in section \ref{initial}.
Then in section \ref{model}, we devise a method to determine the drift and diffusion coefficients for a governing  stochastic differential equation,  with the observation data on the most probable transition trajectory.  To this end, we need a function basis for the drift, and we pick the polynomials basis up to fifth order (or any finite order), together with trigonometric functions. This
  optimization problem is set up in section \ref{opt problem} with an appropriate cost with penalty.
 We devise an algorithm to recursively estimate the coefficients in the diffusion model by feat of the covnex optimization package cvx in Matlab, and then sparse the estimators by vanishing the ones with absolute value smaller than a threshold, in section \ref{algorithm}.
In section \ref{eg}, we demonstrate our method with numerical experiments, and conclude with discussions in section \ref{discussion}.

\section{The Onsager-Machlup function and the most probable transition trajectory}\label{initial}
We start from a scalar stochastic differential equation (SDE)  \cite{Duan}   as the model for the governing stochastic law
\begin{eqnarray}
dx(t) &=& f(x(t)) dt + \varepsilon dB_t,
\end{eqnarray}
where $B_t$ is a Brownian motion defined on a sample space $\Omega$ with probability $\mathbb{P}$,  and $\varepsilon$ is the  noise intensity (i.e., diffusion coefficient).   We are interested in the transition phenomenon between two metastable states, i.e.,  we  examine  random sample trajectories  satisfying the two-point boundary conditions
$$
  x(0) = x_0 \in \mathbb{R}, \; \quad x(T) = x_T \in \mathbb{R},
$$
with a   finite transition time $T$.

The most probable transition  trajectory, denoted by $z(t)$, for this stochastic differential equation,  from state $x_0$ to $x_T$ within a time period $[0, T]$, is the pathway which minimizes  the  Onsager-Machlup (OM)  action functional,  $\int_0^T OM(\dot z(t), z(t)) dt$.   Here the integrand (like a Lagrangian function in classical mechanics)  $OM(\dot z(t), z(t)) $ is known as \cite{Durr, Onsager, Tisza}:
\begin{eqnarray}
OM(\dot z, z) &=& \Big(\frac{f(z) - \dot z}{\varepsilon}\Big)^2  + f'(z).
\end{eqnarray}
We recall the notions   $\dot z = \frac{dz}{dt}$ and $f'(z) = \frac{df}{dz}$.
By the variational principle, we obtain the Euler-Lagrange equation \cite{Durr} to be satisfied by the most probable transition  trajectory $z(t)$:
\begin{eqnarray}
\frac{d}{dt}\frac{\partial OM(\dot z, z)}{\partial \dot z} &=& \frac{\partial OM(\dot z, z)}{\partial z}.
\end{eqnarray}
Thus, we get a second order ODE for the most probable transition trajectory  $z(t)$,  with two-point boundary conditions
\begin{eqnarray}\label{2ndODE}
\ddot z(t) &=& \frac{\varepsilon^2}{2}f''(z(t)) + f'(z(t)) f(z(t)), \quad\quad
z(0) = x_0, \;\; z(T) = x_T.
\end{eqnarray}

As discussed in the previous section, the most probable transition trajectory $ z$  may be observed      in certain dynamical systems. With the observed most probable transition trajectory, we will extract the underlying governing stochastic model in the next section. That is, we will determine the drift $f$ and the diffusion coefficient $\varepsilon$.

\section{Sparse identification of the governing  stochastic differential equation}

\subsection{The equation for the most probable transition trajectory}\label{model}
Based on a time-series data collected with time step  $\delta t$, that is,  $z^{ob}(t), t = 0, \delta t, 2\delta t, \cdots, T$, we want to identify the drift term $f(x)$ and the diffusion coefficient $\varepsilon$.

 We construct a basis library $\vec\Theta(X)$ consisting of polynomial and trigonometric function of observational feature, i.e. $\vec\Theta(X) = \{1, x, x^2, \cdots, x^p, \sin{x}, \cos{x}, \cdots, \sin{qx}, \cos{qx}\}$.  To illustrate our method, we take $p = 5$ and $q=2$. The drift term can then be expressed  as
\begin{eqnarray}
f(x) = \beta_0 + \beta_1 x + \beta_2 x^2 + \beta_3 x^3 + \beta_4 x^4 + \beta_5 x^5 + \beta_6\sin{x} +\beta_7\cos{x}+\beta_8\sin{2x}+\beta_9\cos{2x}. \nonumber\\ \label{f}
\end{eqnarray}
Putting this expression into \eqref{2ndODE},  we then have  the second-order ODE
\begin{eqnarray}\label{2ndODEexpansion}
\ddot z &=& b_0 + b_1 z + b_2 z^2 + b_3 z^3 + \cdots + b_{9} z^{9}+ \sum\limits_{j=1}^4 \big[b_{js}\sin(jz)+b_{jc}\cos(jz)\big] \nonumber\\&&+ \sum\limits_{k=1}^{5} \big[b_{1k} x^k \sin(z) +  b_{2k} x^k \cos(z)
+  b_{3k} x^k \sin(2z) +  b_{4k} x^k \cos(2z)\big],
\end{eqnarray}
with coefficients
\begin{eqnarray}\label{b_beta1}
b_0 &=& \beta_0\beta_1 +  \varepsilon^2 \beta_2,  \quad\quad\quad\quad\quad\quad\quad\quad\quad\quad\quad\;\;
b_5 = 6\beta_1\beta_5 + 6\beta_2\beta_4 +  3\beta_3^2,  \nonumber\\
b_1 &=& 2\beta_0\beta_2 + \beta_1^2 + 3 \varepsilon^2 \beta_3,  \quad\quad\quad\quad\quad\quad\quad\quad\;
b_6 = 7\beta_2\beta_5 + 7\beta_3\beta_4,  \nonumber\\
b_2 &=&  3\beta_0\beta_3 + 3\beta_1\beta_2 + 6 \varepsilon^2 \beta_4, \quad\quad\quad\quad\quad\quad\quad
b_7 = 8\beta_3\beta_5 + 4 \beta_4^2,  \quad\quad\quad\quad\;\; \nonumber\\
b_3 &=& 4\beta_0\beta_4 + 4\beta_1\beta_3  + 2 \beta_2^2 +10\varepsilon^2\beta_5,\quad\quad\quad\;\;\,
b_8 =  9\beta_4\beta_5,  \nonumber\\
b_4 &=&  5\beta_0\beta_5 + 5\beta_1\beta_4 + 5 \beta_2 \beta_3, \quad\quad\quad\quad\quad\quad\quad
b_9 = 5\beta_5^2, \nonumber\\
\end{eqnarray}
\begin{eqnarray}\label{b_beta2}
b_{1s} &=& -\varepsilon^2\beta_6/2 - \beta_0\beta_7 + \beta_1\beta_6-\beta_6\beta_8/2-\beta_7\beta_9/2, \quad\quad
b_{3s} = 3(\beta_6\beta_8-\beta_7\beta_9)/2, \nonumber\\
 b_{1c} &=& -\varepsilon^2\beta_7/2 - \beta_0\beta_6 + \beta_1\beta_7+\beta_7\beta_8/2-\beta_6\beta_9/2, \quad\quad
 b_{3c} = 3(\beta_7\beta_8+\beta_6\beta_9)/2,\nonumber \\
 b_{2s} &=& -2\varepsilon^2\beta_8 - 2\beta_0\beta_9 + \beta_1\beta_8+\beta_6^2/2-\beta_7^2/2, \quad\quad\quad\;\;\;
 b_{4s} = \beta_8^2-\beta_9^2,\nonumber\\
 b_{2c} &=& -2\varepsilon^2\beta_9 + 2\beta_0\beta_8 + \beta_1\beta_9+\beta_6\beta_7,\quad\quad\quad\quad\quad\quad\quad\;
 b_{4c} = 2\beta_8\beta_9.
\end{eqnarray}
and for $i = 1, 2, 3, 4$,
 \begin{eqnarray}\label{b_beta3}
 b_{1i} &=& (i+1)\beta_{i+1}\beta_6-\beta_i\beta_7, \quad
 b_{2i} = (i+1)\beta_{i+1}\beta_7+\beta_i\beta_6,  \quad b_{15} =-\beta_5\beta_7,\quad  b_{25} = \beta_5\beta_6, \nonumber \\
 b_{3i} &=& (i+1)\beta_{i+1}\beta_8-2\beta_{i}\beta_9, \;\,
 b_{4i} = (i+1)\beta_{i+1}\beta_9+2\beta_i\beta_8,\;\, b_{35} =-2\beta_5\beta_9, \;\,
 b_{45} = 2\beta_5\beta_8. \nonumber \\
 \end{eqnarray}
Let $^T$ means transposition, we denote the vectors
\begin{eqnarray*}
\vec{b} &=& (b_0, b_1, b_2, \cdots, b_9, b_{1s}, b_{1c}, b_{2s}, b_{2c}, b_{3s}, b_{3c}, b_{4s}, b_{4c}, b_{11}, b_{21}, b_{31}, b_{41}, \cdots, b_{15}, b_{25}, b_{35}, b_{45})^T, \\
\vec{\beta} &=& (\beta_0, \beta_1, \beta_2, \beta_3, \beta_4, \beta_5, \beta_6, \beta_7, \beta_8, \beta_9)^T.
 \end{eqnarray*}
 For $i = 0, \cdots, 9, 1s, 1c, 2s, 2c, 3s, 3c, 4s, 4c, 11, 21, 31, 41, \cdots, 15, 25, 35, 45$, denote the r.h.s of the equalities in \eqref{b_beta1}-\eqref{b_beta3} as $g_i(\vec{\beta}, \varepsilon)$, and moreover \[\vec{G}(\vec{\beta}, \varepsilon) = (g_0,  \cdots, g_{9}, g_{1s}, g_{1c}, \cdots g_{4s}, g_{4c}, g_{11}, g_{21}, g_{31}, g_{41}, \cdots, g_{51}, g_{51}, g_{51}, g_{51})^T.\] Then the coefficient vector  $\vec{\beta}$ in \eqref{f} is nonlinear related to the coefficients $\vec{b}$ in \eqref{2ndODEexpansion}, and we indicate  this fact by
\begin{equation}\label{b=g(beta)}
\vec{b} = \vec{G}(\vec{\beta}, \varepsilon).
\end{equation}

\subsection{The optimization problem}\label{opt problem}
We use simulated data for the most probable transition trajectory $z(t)$. That is,  we  collect a time history of observation with time step $\delta t$, and denote $z_i^{ob} = z^{ob}(t_{i})$, with $i = 1, 2, \cdots, N+1$, $N = T/\delta t$ and $t_i = (i-1) \delta t$.
Taking center-difference to approximate the second-order accurate solution of $\ddot z$, we get a series of
$\ddot z^{ob}_i = (z^{ob}_{i+1} -2 z^{ob}_i + z^{ob}_{i-1})/\delta t^2$, with $i = 2, \cdots, N$, and $\ddot z^{ob}_1 = (2z^{ob}_1 -5z^{ob}_2 +4z^{ob}_3 -z^{ob}_4)/\delta t^2$, $\ddot z^{ob}_{N+1} = (2z^{ob}_{N+1} -5z^{ob}_N +4z^{ob}_{N-1} -z^{ob}_{N-2})/\delta t^2$. Each pair of data $(z_i^{ob}, \ddot z_i^{ob}), i = 1, 2, 3, \cdots, N+1$, satisfies the ODE \eqref{2ndODEexpansion}, that is
\begin{eqnarray}\label{1stODE_i}
\ddot z_i^{ob} &=& b_0 + b_1 z_i^{ob} + b_2 (z_i^{ob})^2 + \cdots + b_9 (z_i^{ob})^{9} + \sum\limits_{j=1}^4\big[b_{js} \sin(jz_i^{ob})+ b_{jc} \cos(jz_i^{ob})\big] \nonumber\\ &&+\sum\limits_{k=1}^{5} \Big[b_{1k} (z_i^{ob})^k \sin(z_i^{ob}) + b_{2k} (z_i^{ob})^k \cos(z_i^{ob})
+  b_{3k} (z_i^{ob})^k \sin(2z_i^{ob}) +  b_{4k} (z_i^{ob})^k \cos(2z_i^{ob})\Big].\nonumber\\
\end{eqnarray}
Define
\begin{equation}\label{Z-Theta}
\vec{\ddot Z} =
\left(
\begin{array}{c}
           \ddot z_1^{ob} \\
           \ddot z_2^{ob} \\
           \vdots \\
           \ddot z_{N+1}^{ob} \\
\end{array}
\right)
  \mbox{ and }
\vec{\Theta}  =
\left(
\begin{array}{ccccccc}
1 & z_1^{ob} & \cdots & (z_1^{ob})^9 & \mathbf{B}_{0}(z_1^{ob})& \cdots & \mathbf{B}_{5}(z_1^{ob})\\
 1 & z_2^{ob} &  \cdots & (z_2^{ob})^9 & \mathbf{B}_{0}(z_2^{ob})& \cdots & \mathbf{B}_{5}(z_2^{ob})\\
 \vdots & \vdots & \ddots & \vdots & \vdots & \ddots & \vdots \\
 1 & z_{N+1}^{ob} &  \cdots & (z_{N+1}^{ob})^9 &\mathbf{B}_{0}(z_{N+1}^{ob})& \cdots & \mathbf{B}_{5}(z_{N+1}^{ob}))
 \end{array}
 \right),
\end{equation}
with $\mathbf{B}_{0}(z) = [\sin(z), \cos(z), \sin(2z), \cos(2z), \sin(3z), \cos(3z), \sin(4z), \cos(4z)]$ and $\mathbf{B}_{k}(z) = [z^{k}\sin(z), z^{k}\cos(z), z^{k}\sin(2z), z^{k}\cos(2z)]$, $k = 1, 2, 3, 4, 5$.
The matrix formulation for the `unknown' $b$  is
 \begin{eqnarray}\label{matrix}
  \vec{\ddot Z}  &=& \vec{\Theta} \vec{b}.
 \end{eqnarray}
This is usually an overdetermined system consisting of $N+1$ equalities and 38 unknown coefficients. Remember that, our ultimate aim is extracting a sparse $\vec{\beta}$ and a $\varepsilon$ from
\begin{eqnarray}\label{matrix_arm}
  \vec{\ddot Z}  &=& \vec{\Theta}\vec{G}(\vec{\beta}, \varepsilon).
 \end{eqnarray}
 Starting from the most probable trajectory satisfying \eqref{2ndODE}, we initially estimate the coefficients $\vec{b}$ in \eqref{matrix}, and then work out $\vec{\beta}$ and $\varepsilon$ from the relationship \eqref{b=g(beta)}.
The sparse problem is
\begin{eqnarray}\label{Opt}
&&\min ||\vec{\beta}||_0 \nonumber \\
&&\mathrm{s. t.}\;\;
\begin{cases}
 \vec{\ddot Z}  = \vec{\Theta} \vec{b},\\
  \vec{b} = \vec{G}(\vec{\beta}, \varepsilon).
  \end{cases}
 \end{eqnarray}
For $\vec{x} = (x_1, \cdots, x_n)$, we denote the $L_2$ norm $||\vec{x}||^2_2 = \sum_{i = 1}^{i = n} x_i^2$,  the $L_1$ norm $||\vec{x}||_1= \sum_{i = 1}^{i = n} |x_i|$ and the $L_0$ norm $||\vec{x}||_0$ the number of non-zero $x_i$. The sparse problem can be transformed into the following regularized least-square problem with an $L_0$ penalty \cite{Xu2012}
\begin{eqnarray}\label{LS+R0}
\min\limits_{\vec{\beta}, \varepsilon, \vec{b}} \; || \vec{\ddot Z}  - \vec{\Theta} \vec{b}||^2_2 + \kappa_1||\vec{b} - \vec{G}(\vec{\beta}, \varepsilon) ||^2_2 + \kappa_2||\vec{\beta}||_0,
 \end{eqnarray}
with positive $\kappa_1, \kappa_2$. The first two terms in this optimization problem are a weighted sum of the inconsistency of the two constraint functions in \eqref{Opt}. Note that, we cannot estimate $\vec{\beta}$ from \eqref{LS+R0} directly due to nonlinearity in $\vec{G}(\vec{\beta}, \varepsilon)$.
 We introduce an intermediate variable $\vec{b}$, which is then linear in the constraint equations in \eqref{Opt}, and we intend to estimate it firstly. As sparse $\vec{\beta}$ results in sparse $\vec{b}$, we substitute $||\vec{\beta}||_0$ in \eqref{LS+R0} to $||\vec{b}||_0$ for sparsing $\vec{b}$. To make the optimization problem for $\vec{b}$ solvable, we furthermore relax $||\vec{b}||_0$ to $||\vec{b}||_2$, and take
\begin{eqnarray}\label{LS+R2}
 || \vec{\ddot Z}  - \vec{\Theta} \vec{b}||^2_2 + \kappa_1||\vec{b} - \vec{G}(\vec{\beta}, \varepsilon) ||^2_2 + \kappa_2||\vec{b}||_2^2
 \end{eqnarray}
as loss function for $\vec{b}$.

 \subsection{An implementable algorithm for the optimization problem}\label{algorithm}


  We start from the linear system $\vec{\Theta} \vec{b} = \vec{\ddot Z}$. The analytical solution $\vec{b} = (\vec{\Theta}^T\vec{\Theta})^{-1}(\vec{\Theta}^T\vec{\ddot Z})$ relies on the reversibility of $\vec{\Theta}^T\vec{\Theta}$. By the help of the backslash, $\vec{b} = \vec{\Theta}\setminus \vec{\ddot Z}$ is more prone to result in sparse solution in an overdetermined system\cite{Kutz}. In our experiment, we found that the above two kinds of solutions are usually associated with the rank deficiency, so we employ a convex optimization package, cvx, to get an initial estimator $\vec{b}_{init}$. We then solve $\vec{\beta}_{init}$ and $\varepsilon_{init}$ from \eqref{b=g(beta)}.

  We then take a loop to update $\vec{b}$ in the next improvement step, and solve $\vec{\beta}$ and $\varepsilon$, until the newest two sets of solutions are close enough or the loop time reaches a setted number. Here $\vec{b}$ is updated by the cvx package for
  $\min_{\vec{b}} \; || \vec{\ddot Z}  - \vec{\Theta} \vec{b}||^2_2 + \kappa_1||\vec{b} - \vec{G}(\vec{\tilde\beta}, \tilde\varepsilon) ||^2_2 + \kappa_2||\vec{b}||_2^2$, with $\vec{\tilde\beta}$ and $\tilde\varepsilon$ in $\vec{G}(\vec{\tilde\beta}, \tilde\varepsilon)$ the newly solved solutions in the last loop.


 For each estimator $\vec{\beta}$ solved from the updated $\vec{b}$, we vanish $\vec{\beta}$'s elements whose absolute value is smaller than the threshold $\theta_T$ as in \cite{Brunton, Rudy}. When the adjacent two estimators are close enough, we terminate the update loop, calculate errors $E1, E2, \cdots, E6$, which will be given below, compare $E6$ with the old one and generate a new $\theta_T$ according its renewal rule. The best $\theta_T$ is selected related to the smallest $E3$. The weight parameters $\kappa_1$ and $\kappa_2$ are chosen by cross validation, trading off of errors $E1$ and $E3$, and considering of the convergence of the estimators.

 We first give a drift function $f(x)$ and a diffuse coefficient $\varepsilon$. Use them to simulate a time-series date $z^{ob}_i, i=1, 2, \cdots, N+1$ from the Euler-Lagrange equation of the SDE by a  shooting method. We then approximate the time-series date $\ddot z$ from $z^{ob}_i$ by center-difference, and get $\mathbf{\ddot Z} = [\ddot z^{ob}_1, \ddot z^{ob}_2, \cdots, \ddot z^{ob}_{N+1}]^T$.
 Second, we construct a basis library for Euler-Lagrange equation, the following basis corresponds to a basis expression of $f$ as in \eqref{f}.
  \begin{eqnarray*}
  \vec{\Theta}(z) &=& \{1, z, z^2, \cdots, z^{9}, \sin{z}, \cos{z}, \sin{2z}, \cos{2z}, \sin{3z}, \cos{3z}, \sin{4z}, \cos{4z},\\
   &&  z\sin{z}, z\cos{z}, z\sin{2z}, z\cos{2z}, \cdots, z^5\sin{z}, z^5\cos{z}, z^5\sin{2z}, z^5\cos{2z}\},
  \end{eqnarray*}
   and develop the basis to get a matrix $\vec{\Theta}$ of time history as in \eqref{Z-Theta}.
   Thirdly, we split the data randomly to training part (70\% data) and testing part (the left 30\%). Accordingly, $\vec{\ddot Z}$ and $\vec{\Theta}$ are divided into  $\vec{\ddot Z}_{train}$, $\vec{\ddot Z}_{test}$ and $\vec{\Theta}_{train}$, $\vec{\Theta}_{test}$, respectively.

   After the preparations, we initially get initial estimators and errors.
   \begin{algorithm}
   \caption{ Get initial estimators and errors}
   \begin{algorithmic}
   \STATE{1. get $\vec{b}_{init}=\min\limits_{\vec{b}}{||\vec{\ddot Z}_{train}-\vec{\Theta}_{train} \vec{b}||^2_2}$ by cvx,}
   \STATE{2. solve out $\vec{\beta}_{init}$ and $\varepsilon_{init}$ by findbeta1.m in Algorithm \ref{findbeta1},}
   \STATE{3. calculate $\vec{G0}=\vec{G}(\vec{\beta}_{init}, \varepsilon_{init})$,}
   \STATE{4. calculate the test errors \\
  $E1_{init} = ||\vec{\Theta}_{test} \vec{b}_{init} - \vec{\ddot Z}_{test} ||_2^2$, \quad\quad\quad\quad\quad  $E2_{init}=|| \vec{b}_{init} - \vec{G}(\vec{\beta}_{init}, \varepsilon_{init})||_2^2$,\\
   $E3_{init} = ||\vec{\Theta}_{test} \vec{G}(\vec{\beta}_{init} \varepsilon_{init}) - \vec{\ddot Z}_{test} ||_2^2$, \quad $E4_{init}=||\vec{\beta}_{init}||_0$, \quad $E5_{init}=||\vec{b}_{init}||_0$.}
   \end{algorithmic}
   \end{algorithm}

  \begin{breakablealgorithm}
   \caption{ The main loop}
   \begin{algorithmic}
   \STATE{Set a value set $\Omega_1$ for $\kappa_1$, and $\Omega_2$ for $\kappa_2$, $BetaTotal=\mathrm{cell}(1, 1)$, $EBest=\mathrm{cell}(1, 1)$;}
   \FOR {$k1=1:length(\Omega_1)$}
   \FOR {$k2=1:length(\Omega_2)$}
   \STATE $\kappa_1=\Omega_1(k1), \kappa_2=\Omega_2(k2);$
   \STATE $\vec{b}_{best} = \vec{b}_{init}$, $\vec{\beta}_{best} = \vec{\beta}_{init}$, $\varepsilon_{best} = \varepsilon_{init}$,\\
  $E1_{best}=E1_{init}$, $E2_{best}=E2_{init}$, $E3_{best}=E3_{init}$, $E4_{best}=E4_{init}$, $E5_{best}=E5_{init}$, \\$E6_{best}=E1_{best}+\kappa_1 E2_{best} +\kappa_2 E5_{best}$,\\
  $Ebest=[E1_{best}, E2_{best}, E3_{best}, E4_{best}, E5_{best}, E6_{best}, 0, \kappa_1, \kappa_2];$
   \STATE set an initial $\theta_T$ and its initial step $h$
   \FOR {Iter=1:TT}
   \STATE $\vec{Gb}=\vec{G0}$, $\vec{\beta_0}=\vec{\beta}_{init}$,
    $BetaTotal\{k1, k2\}\{Iter, 1\}=[\vec{\beta_0}; \vec{\varepsilon}]$; \\
    $BetaTotal\{k1, k2\}\{Iter, 2\}=E1_{best}$,  $BetaTotal\{k1, k2\}\{Iter, 3\}=E2_{best}$,  \\
    $BetaTotal\{k1, k2\}\{Iter, 4\}=E3_{best}$,  $BetaTotal\{k1, k2\}\{Iter, 5\}=E4_{best}$,  \\
    $BetaTotal\{k1, k2\}\{Iter, 6\}=E5_{best}$,  $BetaTotal\{k1, k2\}\{Iter, 7\}=E6_{best}$,  \\
    $BetaTotal\{k1, k2\}\{Iter, 8\}=0$,  $BetaTotal\{k1, k2\}\{Iter, 9\}=\kappa_1$, \\
    $BetaTotal\{k1, k2\}\{Iter, 10\}=\kappa_2$; \\
   \FOR{Iterations=1:T}
   \STATE $non0Gb=(\mathrm{abs}(\vec{Gb})>0), numb=\mathrm{sum}(non0Gb), b=\mathrm{zeros}(38, 1)$; \\
    Get $numb$ dimensional $\vec{bu}$ minimizing\\ ${||\vec{\ddot Z}_{train}-\vec{\Theta}_{train}(:, non0Gb)\vec{bu}||^2_2}+ \kappa_1||\vec{bu} - \vec{Gb}(non0G, 1) ||^2_2 + \kappa_2||\vec{bu}||_2^2$ by cvx;\\
    $b(non0Gb, 1)=bu;$\\
    get $\vec{\tilde\beta}$ and $\tilde\varepsilon$ from $\vec{b}$ by findbeta1.m, $small=(\mathrm{abs}(\vec{\tilde\beta})<\theta_T), \vec{\tilde\beta}(small, 1)=0;$\\
    calculate $\vec{Gb}=\vec{G}(\vec{\tilde\beta}, \tilde\varepsilon);$\\
    get  $\vec{\beta}$ and $\varepsilon$ by findbeta0.m; \\
    calculate $\vec{Gb}=\vec{G}(\vec{\beta}, \varepsilon)$;\\
    \IF{$\text{norm}{(\vec{\beta}-\vec{\beta_0}, 1)}<\epsilon_1$}
    \STATE  break
    \ENDIF
    \STATE $\vec{\beta_0}=\vec{\beta}$;
   \ENDFOR
   \STATE  $E1 = ||\vec{\Theta}_{test} \vec{b} - \vec{\ddot Z}_{test} ||_2^2,$ \quad  $E2=|| \vec{b} - \vec{G}(\vec{\beta}, \varepsilon)||_2^2$,\quad
   $E3 = ||\vec{\Theta}_{test} \vec{G}(\vec{\beta}, \varepsilon) - \vec{\ddot Z}_{test} ||_2^2,$ \\
   $E4=||\vec{\beta}||_0$, \quad\quad\quad\quad\quad\quad\, $E5=||\vec{b}||_0$, \quad\quad\quad\quad\quad $E6=E1+\kappa_1 E2 +\kappa_2 E5$, \\
   $BetaTotal\{k1, 1\}\{Iter, 1\}=[\vec{\beta}; \varepsilon]$,
   $BetaTotal\{k1, k2\}\{Iter, 2\}=E1$, \\ $BetaTotal\{k1, k2\}\{Iter, 3\}=E2$,
    $BetaTotal\{k1, k2\}\{Iter, 4\}=E3$, \\ $BetaTotal\{k1, k2\}\{Iter, 5\}=E4$,
    $BetaTotal\{k1, k2\}\{Iter, 6\}=E5$, \\ $BetaTotal\{k1, k2\}\{Iter, 7\}=E6$,
    $BetaTotal\{k1, k2\}\{Iter, 8\}=\theta_T$, \\ $BetaTotal\{k1, k2\}\{Iter, 9\}=\kappa_1$,
    $BetaTotal\{k1, k2\}\{Iter, 10\}=\kappa_2$;
   \IF{$E6<E6_{best}$}
   \STATE $E1_{best}=E1$, $E2_{best}=E2$, $E3_{best}=E3$, $E4_{best}=E4$, $E5_{best}=E5$, \\ $E6_{best}=E6$, $E7_{best}=E7$, $\theta_T=\theta_T+h$;\\ $Ebest=[Ebest; E1_{best}, E2_{best}, E3_{best}, E4_{best}, E5_{best}, E6_{best}, \theta_T, \kappa_1, \kappa_2];$
   \ELSE
   \STATE $\theta_T=\max{(0, \theta_T-2h)},\; h=h/2,\; \theta_T=\theta_T+h;$
   \ENDIF
   \IF{$h<\epsilon_2$}
   \STATE break
   \ENDIF
    \ENDFOR
   \STATE $[val3, pos3]=\min{(Ebest(:, 3))}, pos=pos3(1)$,\\ $ET1(k1, k2)=Ebest(pos, 1)$,  $ET2(k1, k2)=Ebest(pos, 2)$,\\ $ET3(k1, k2)=Ebest(pos, 3)$, $ET4(k1, k2)=Ebest(pos, 4)$,\\ $ET5(k1, k2)=Ebest(pos, 5)$, $ET6(k1, k2)=Ebest(pos, 6)$, $THR=Ebest(pos, 7);$\\
   $row=\mathrm{find}(\mathrm{cell2mat}(BetaTotal\{k1, k2\}(:, 2))==THR)$; \\
   $BETA(12(k1-1)+1:12(k1-1)+11, k2)=BetaTotal\{k1, k2\}\{row(1), 1\}(:, \mathrm{end})$; \\
   \ENDFOR
   \ENDFOR
   \end{algorithmic}
   \end{breakablealgorithm}

Recall that our aim is to extract a sparse $\vec{\beta}$ and obtain an estimator $\varepsilon$ from \eqref{matrix_arm}.
Thus, the error $E3=||\ddot Z-\Theta\vec{G}(\vec{\beta}, \varepsilon)||^2_2$ should be small for the estimators $\vec{\beta}$ and $\varepsilon$.
We introduce an intermediate variable $\vec{b}$, which can be estimated initially from \eqref{matrix}, such that the error $E1=||\ddot Z-\Theta\vec{b}||^2_2$ reaches a minimum, and in the following update steps the estimator $\vec{b}$ is got by minimizing \eqref{LS+R2}. In \eqref{LS+R2}, the second term $||\vec{b}-\vec{G}(\vec{\beta}, \varepsilon)||^2_2$  tends to draw $\vec{b}$ toward $\vec{G}(\vec{\beta}, \varepsilon)$, and the regularizer $||\vec{b}||^2_2$ apts to results in some elements of $\vec{b}$ toward to 0. The bigger the weights $\kappa_1$ and $\kappa_2$ are, the stronger the two terms effect. These simultaneously sacrifice the error $E1$. Therefore, for a set of $\kappa_1, \kappa_2, \theta_T$, after update the estimators $\vec{b}, \vec{\beta}$ and $\varepsilon$ by minimizing $E6$, we select the results which makes $E6$ smaller and smaller and put them in the best sets $Ebest$. We take the results corresponding to the smallest $E3$ from $Ebest$.
After the main loop, each $(\kappa_1, \kappa_2)$ corresponding to an estimator and some errors.
We then pick   the results with $E3$ below a certain level, with  trade off  between $E1$ and $E3$, and take a convergence estimator.

Based on the relationships between $\vec{b}$ and $(\vec{\beta}, \varepsilon)$ \eqref{b_beta1}--\eqref{b_beta3}, the computing processes of findbeta1.m is described in Algorithm \ref{findbeta1}:
\begin{breakablealgorithm}
\caption{findbeta1.m}\label{findbeta1}
\begin{algorithmic}
 \STATE{start from $b_9$, get a real number $\beta_5=\sqrt{\bar b_9/5}$, with $\bar b_9=\max{(b_9, 0)}$}
 \IF{$\beta_5\neq 0$}
 \STATE $\beta_4=b_8/(9\beta_5)$, $\beta_3=(b_7-4\beta_4^2)/(8\beta_5)$, $\beta_2=(b_6-7\beta_3\beta_4)/(7\beta_5)$, \\ $\beta_1=(b_5-3\beta_3^2-6\beta_2\beta_4)/(6\beta_5)$,\\
 then the equations $b_i$, $i =0, 1, 2, 3, 4$ are linear with respect to $\varepsilon^2$ and $\beta_0$, \\set
 $A5=[\beta_2\;\beta_1; 3\beta_3\; 2\beta_2; 6\beta_4\; 3\beta_3; 10\beta_5\; 4\beta_4; 0\; 5\beta_5]$,\\
 $R5=[b_0; b_1-\beta_1^2; b_2-3\beta_1*\beta_2; b_3-2\beta_2^2-4\beta_1\beta_3; b_4-5\beta_1\beta_4-5\beta_2\beta_3]$,\\
 then $x5=A5\backslash R5$, $e2=x5(1)$, $\varepsilon=\sqrt{\mathrm{abs}(e2)}$,\\
 \IF{$e2<0$}
 \STATE $\beta_i=-\beta_i$, for $i=5, 4, 3, 2, 1, 0,$
 \ENDIF
 \ELSE 
 \STATE $\beta_4=\sqrt{\bar b_7/4}$, with $\bar b_7=\max{(b_7, 0)},$
 \IF{$\beta_4\neq 0$}
 \STATE $\beta_3=b_6/(7\beta_4),$
        $\beta_2=(b_5-3\beta_3^2)/(6\beta_4),$
        $\beta_1=(b_4-5\beta_2\beta_3)/(5\beta_4),$\\
        $A4=[\beta_2\; \beta_1; 3\beta_3\; 2\beta_2; 6\beta_4\; 3\beta_3; 0\; 4\beta_4],$\\
        $R4=[b_0; b_1-\beta_1^2; b_2-3\beta_1\beta_2; b_3-2\beta_2^2-4\beta_1\beta_3],$\\
        then $x4=A4\backslash R4, \beta_0=x4(2), e2=x4(1), \varepsilon=\sqrt{\mathrm{abs}(e2)},$
 \IF{$e2<0$}
 \STATE $\beta_i=-\beta_i$, for $i=4, 3, 2, 1, 0,$
 \ENDIF
 \ELSE
 \STATE $\beta_3=\sqrt{\bar b_5/3}$, with $\bar b_5=\max{(b_5, 0)},$
 \IF{$\beta_3 \neq 0$}
 \STATE $\beta_2=b_4/(5\beta_3)$, $\beta_1=(b_3-2\beta_2^2)/(4\beta_3)$, \\
 $A3=[\beta_2\; \beta_1; 3\beta_3\; 2\beta_2; 0\; 3\beta_3]$,
            $R3=[b_0; b_1-\beta_1^2; b_2-3\beta_1\beta_2]$,\\
            then $x3=A3\backslash R3, \beta_0=x3(2), e2=x3(1), \varepsilon=\sqrt{\mathrm{abs}(e2)}$;
 \IF{$e2<0$}
 \STATE $\beta_i=-\beta_i$, for $i=3, 2, 1, 0,$
 \ENDIF
 \ELSE
 \STATE $\beta_2=\sqrt{\bar b_3/2}$, with $\bar b_3=\max{(b_3, 0)},$
 \IF{$\beta_2\neq 0$}
 \STATE $\beta_1=b_2/(3\beta_2), \beta_0=(b_1-\beta_1^2)/(2\beta_2), e2=(b_0-\beta_0\beta_1)/\beta_2, \varepsilon=\sqrt{\mathrm{abs}(e2)},$
 \IF{$e2<0$}
 \STATE $\beta_i=-\beta_i$, for $i=2, 1, 0$,
 \ENDIF
 \ELSE
 \STATE $\beta_1=\sqrt{\bar b_1}$, with $\bar b_1=\max{(b_1, 0)},$
 \IF{$\beta_1\neq 0$}
 \STATE $\beta_0=b_0/\beta_1$, and set $\varepsilon=10$,
 \ELSE
 \STATE set $\beta_0=10, \varepsilon=10,$
 \ENDIF
 \ENDIF
 \ENDIF
 \ENDIF
 \ENDIF
 \IF{$\beta_i, i=0, 1, 2, 3, 4, 5$ are all 0}
 \IF{$b_{4c} == 0$}
 \IF{$b_{4s} >0$}
 \STATE $\beta_9=0, \beta_8=\sqrt{b_{4s}}, \beta_6=2b_{3s}/(3\beta_8), \beta_7=2b_{3c}/(3\beta_8),$\\
 $A8=[-\beta_6/2\; -\beta_7; -\beta_7/2\;\beta_6; -2\beta_8
 \;0; 0\;2\beta_8],$\\
 $R8=[b_{1s}+b_{3s}/3-\beta_1\beta_6; b_{1c}-b_{3c}/3-\beta_1\beta_7; b_{2s}-(\beta_6^2-\beta_7^2)/2-\beta_1\beta_8; b_{2c}-\beta_6\beta_7-\beta_1\beta_9],$\\
 then $x8=A8\backslash R8, e2=x8(1),$
 \IF{$e2<0$}
 \STATE $\beta_i=-\beta_i,$ for $i=8, 7, 6$,
 \ENDIF
 \ENDIF
 \IF{$b_{4s}<0$}
 \STATE  $\beta_8=0, \beta_9=\sqrt{-b_{4s}}, \beta_6=2b_{3c}/(3\beta_9), \beta_7=-2b_{3s}/(3\beta_9),$\\
 $A9=[-\beta_6/2\; -\beta_7; -\beta_7/2\;\beta_6; 0\; -2\beta_9; 2-\beta_9\;0],$\\
 $R9=[b_{1s}-b_{3s}/3-\beta_1\beta_6; b_{1c}+b_{3c}/3-\beta_1\beta_7; b_{2s}-(\beta_6^2-\beta_7^2)/2-\beta_1\beta_8; b_{2c}-\beta_6\beta_7-\beta_1\beta_9],$\\
 then $x9=A9\backslash R9, e2=x9(1),$
 \IF{$e2<0$}
 \STATE $\beta_i=-\beta_i,$ for $i=9, 7, 6$,
 \ENDIF
 \ENDIF
 \IF{$b_{4s}==0$}
 \STATE $\beta_8=0, \beta_9=0,$
 \IF{$b_{2c}==0$}
 \IF{$b_{2s}>0$}
 \STATE $\beta_6=\sqrt{2b_{2s}}, \beta_7=0, A67=[-\beta_6/2\; -\beta_7; -\beta_7/2\;\beta_6; -2\beta_8\; -2\beta_9; 2-\beta_9\; 2\beta_8]$,\\
 $R67=[b_{1s}-\beta_1\beta_6; b_{1c}-\beta_1\beta_7; b_{2s}-(\beta_6^2-\beta_7^2)/2-\beta_1\beta_8; b_{2c}-\beta_6\beta_7-\beta_1\beta_9],$\\
 then $x67=A67\backslash R67, e2=x67(1), \beta_6=\mathrm{sign}(e2)\beta_6$,
 \IF{$\beta_0==10$}
 \STATE $\beta_0=\mathrm{sign}(e2)x67(2),$
 \ENDIF
 \IF{$\varepsilon==10$}
 \STATE $\varepsilon=\sqrt{\mathrm{abs}(e2)},$
 \ENDIF
 \ENDIF
 \IF{$b_{2s}<0$}
 \STATE $\beta_7=\sqrt{-2b_{2s}}, \beta_6=0, A76=[-\beta_6/2\; -\beta_7; -\beta_7/2\;\beta_6; -2\beta_8\; -2\beta_9; 2-\beta_9\; 2\beta_8]$,\\
 $R76=[b_{1s}-\beta_1\beta_6; b_{1c}-\beta_1\beta_7; b_{2s}-(\beta_6^2-\beta_7^2)/2-\beta_1\beta_8; b_{2c}-\beta_6\beta_7-\beta_1\beta_9],$\\
 then $x76=A76\backslash R76, e2=x76(1), \beta_7=\mathrm{sign}(e2)\beta_7$,
 \IF{$\beta_0==10$}
 \STATE $\beta_0=\mathrm{sign}(e2)x76(2),$
 \ENDIF
 \IF{$\varepsilon==10$}
 \STATE $\varepsilon=\sqrt{\mathrm{abs}(e2)},$
 \ENDIF
 \ENDIF
 \IF{$b_{2s}==0$}
 \STATE $\beta_6=0, \beta_7=0,$
 \ENDIF
 \ELSE
 \STATE $\beta_7=\sqrt{\sqrt{b_{2s}^2+b_{2c}^2}-b_{2s}}, \beta_6=b_{2c}/\beta_7,$ $e2=(b_{1s\beta_6+b_{1c}\beta_7})/(-(\beta_6^2+\beta_7^2)/2)+2\beta_1$, $\beta_6=\mathrm{sign}(e2)\beta_6, \beta_7=\mathrm{sign}(e2)\beta_7$,
 \IF{$\beta_0==10$}
 \STATE $\beta_0=(b_{1c}\beta_6-b_{1s}\beta_7)/(\beta_6^2+\beta_7^2),$
 \ENDIF
 \IF{$\varepsilon==10$}
 \STATE $\varepsilon=\sqrt{\mathrm{abs}(e2)},$
 \ENDIF
 \ENDIF
 \ENDIF
\ELSE
\STATE $\beta_9=\sqrt{\sqrt{b_{4s}^2+b_{4c}^2}-b_{4s}/2}$, $\beta_8=b_{4c}/(2\beta_9)$, $\beta_6=2/3(b_{3c}\beta_9+b_{3s}\beta_8)/(\beta_8^2+\beta_9^2)$, \\ $\beta_7=2/3(b_{3c}\beta_8-b_{3s}\beta_9)/(\beta_8^2+\beta_9^2)$,
$e2=-(b_{2c}-\beta_6\beta_7)/(2\beta_9)$,
\IF{$e2<0$}
\STATE $\beta_i=-\beta_i$, for $i=9, 8, 7, 6,$
\ENDIF
\ENDIF
\ELSE
\STATE $AA=[2\beta_2\; -\beta_1 \;\;0\;\; 0; \beta_1\;\; 2\beta_2\;\; 0\;\; 0; 0\;\; 0\;\; 2\beta_2\; -2\beta_1; 0\;\; 0\;\; 2\beta_1\;\; 2\beta_2; 3\beta_3\; -\beta_2\;\; 0\;\; 0; $ \\ \quad\quad\quad $ \beta_2\;\; 3\beta_3\;\; 0\;\; 0; 0\;\; 0\;\; 3\beta_3\; -2\beta_2; 0\;\; 0\;\; 2\beta_2\;\; 3\beta_3; 4\beta_4\; -\beta_3\;\; 0\;\; 0; \beta_3\;\; 4\beta_4\;\; 0\;\; 0;$\\ \quad\quad\quad $0\;\; 0\;\; 4\beta_4\; -2\beta_3; 0\;\; 0\;\; 2\beta_3\;\; 4\beta_4;5\beta_5\; -\beta_4\;\; 0\;\; 0; \beta_4\;\; 5\beta_5\;\; 0\;\; 0; 0\;\; 0\;\; 5\beta_5\; -2\beta_4; $ \\ \quad\quad\quad $ 0\;\; 0\;\; 2\beta_4\;\; 5\beta_5;  0\; -\beta_5\;\; 0\;\; 0; \beta_5\;\; 0 \;\;0 \;\;0; 0 \;\;0 \;\;0 \;-2\beta_5; 0\;\; 0\;\; 2\beta_5 \;\;0]$, \\
    $RR = [b_{11};b_{21};b_{31};b_{41};b_{12};b_{22};b_{32};b_{42};b_{13};b_{23};b_{33};b_{43};b_{14};b_{24};b_{34};b_{44};b_{15};b_{25};b_{35};b_{45}]$,
    \\ $x69=AA\backslash RR$,
    $\beta_6=x69(1), \beta_7=x69(2), \beta_8=x69(3), \beta_9=x69(4),$
\ENDIF
\end{algorithmic}
\end{breakablealgorithm}

Most procedures in findbeta0.m are the same as in findbeta1.m, besides that the newest estimator $\vec{\beta}$ is used in the next solving process for $\vec{\beta}$ from a more newer updated $\vec{b}$ to maintain the 0 in $\vec{\beta}$.

\section{Numerical experiments}\label{eg}

\begin{description}
  \item[Example] We reconstruct a dynamical system fluctuated by white noise with fixed boundary values
\end{description}
\begin{equation}\label{Ex1}
  d X_t = f(X_t)\ dt + \varepsilon\ dB_t, \quad \quad  X(0) = X_0, \quad X(1) = X_T,
\end{equation}
by a time-series observations of its most probable transition trajectory.
We construct a basis library consisting of polynomial and trigonometric functions of observational feature, $\{1, X, X^2, \cdots, X^5, \sin{X}, \cos{X}, \sin{2X}, \cos{2X}\}$. Then the drift $f(X)$ can be written as in \eqref{f}. The most probable transition trajectory, $z(t)$, satisfies a second order ODE as in \eqref{2ndODE}. Substituting the expression of $f$ into \eqref{2ndODE}, and combining like terms, the second ODE would turn into the form as in \eqref{2ndODEexpansion}--\eqref{b_beta3}.

We reconstruct stochastical systems in three cases:
\begin{enumerate}
\item[I.] $f(X) = 0.5X-1.2X^3 + \sin{X}$, and $\varepsilon=0.8$;
\item[II.] $f(X) = 0.5X-1.2X^3$, and $\varepsilon=0.8$;
\item[III.] $f(X) = \cos{X}$, and $\varepsilon=0.8$.
\end{enumerate}
Taking the boundary values $a_0=0, a_T=\sqrt{2}$, we simulate by a shooting method to get a time-series observation of the most probable trajectory, denoted as $z^{ob}_i = z^{ob}(t_i)$, with $t_i = (i-1) \delta t$, $\delta t = 1/N$, and $i = 1, 2, \cdots, 1001$. Based on the time-series data, we reversely reconstruct the system \eqref{Ex1}, namely, distill the nonzero coefficients in $f$ and estimate them together with the diffusion coefficient $\varepsilon$.

\begin{figure}[htbp]
\begin{center}
\includegraphics[width=5.3cm, height=6cm]{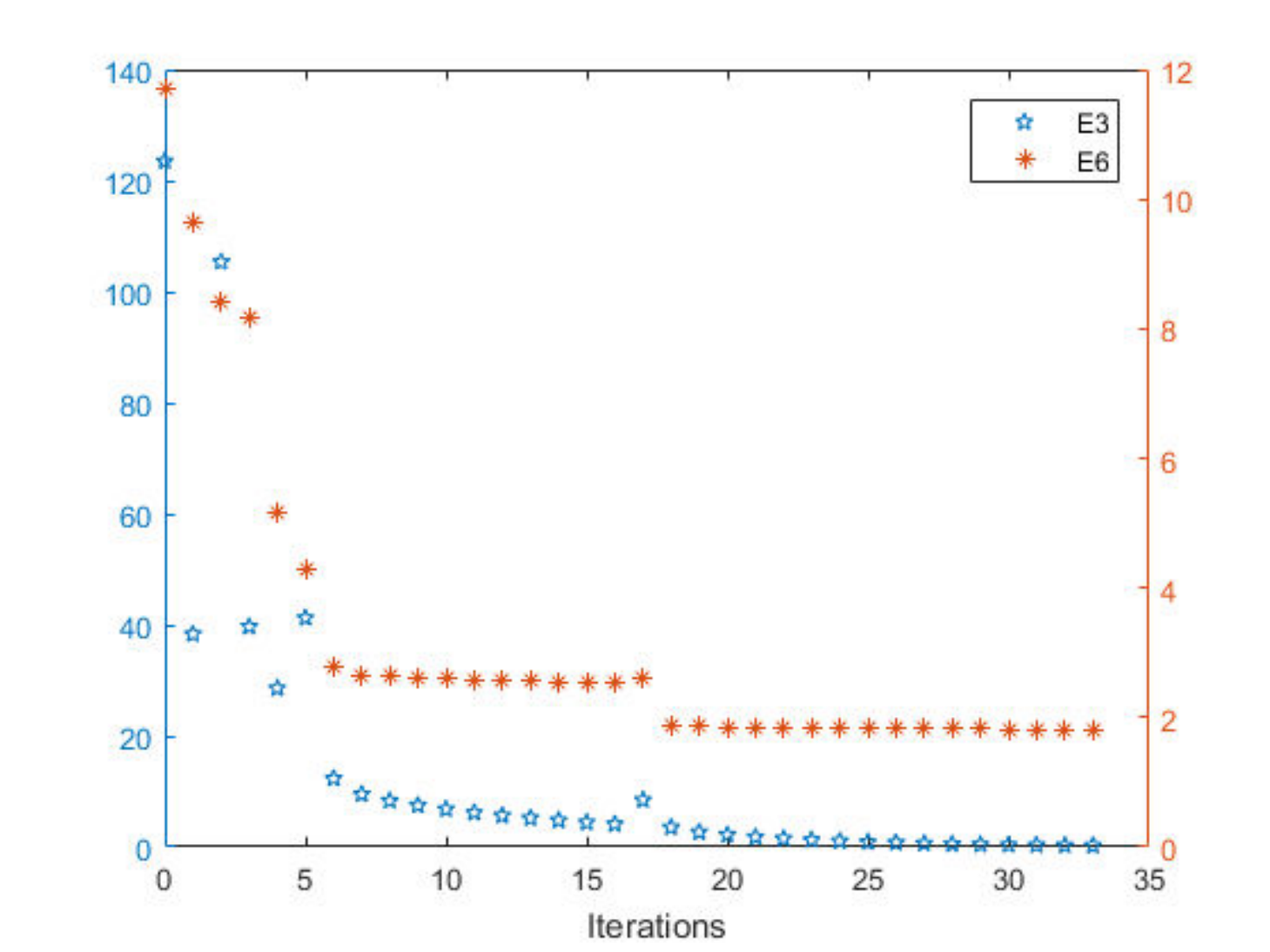}
\includegraphics[width=5.3cm, height=6cm]{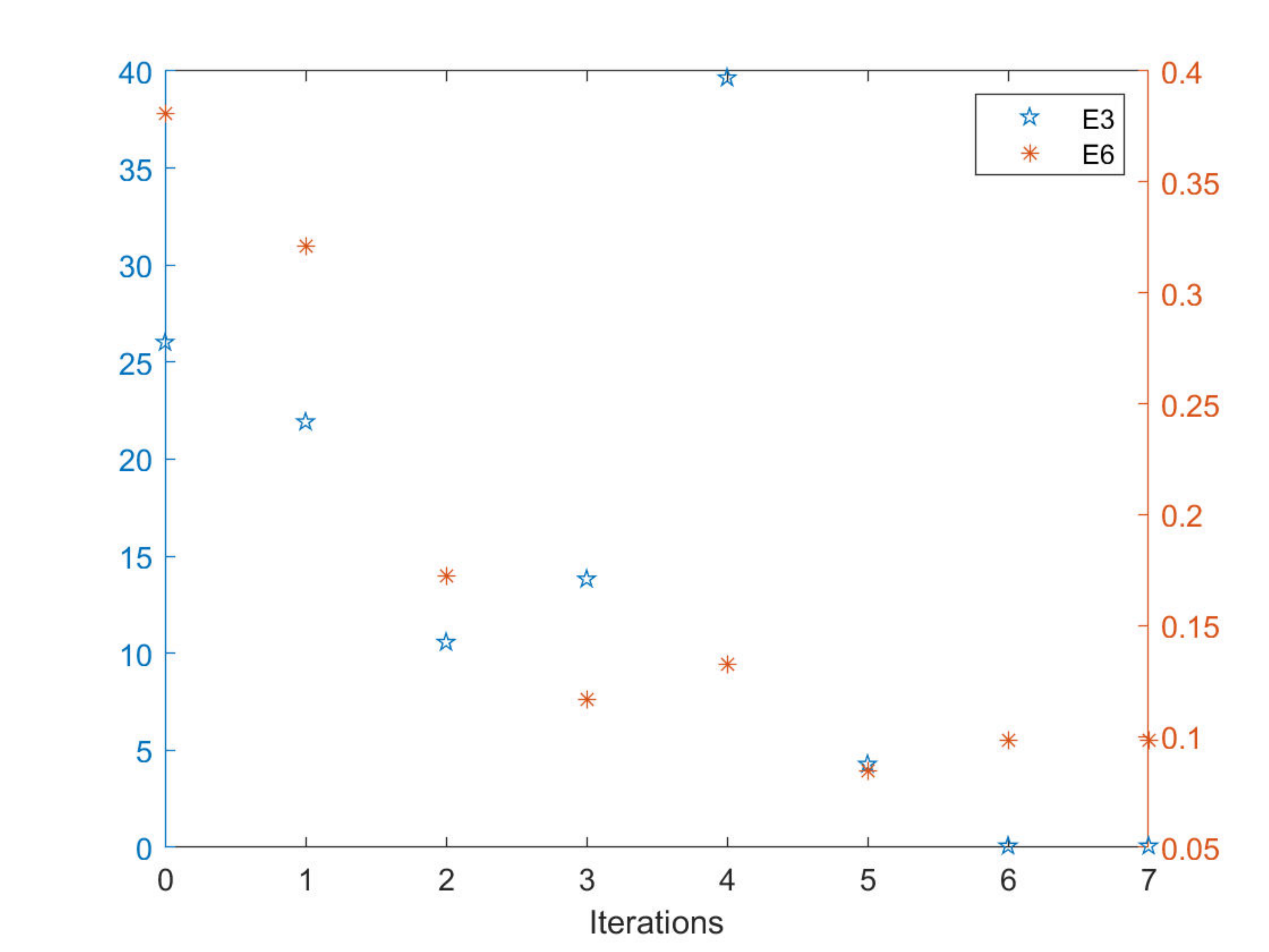}
\includegraphics[width=5.3cm, height=6cm]{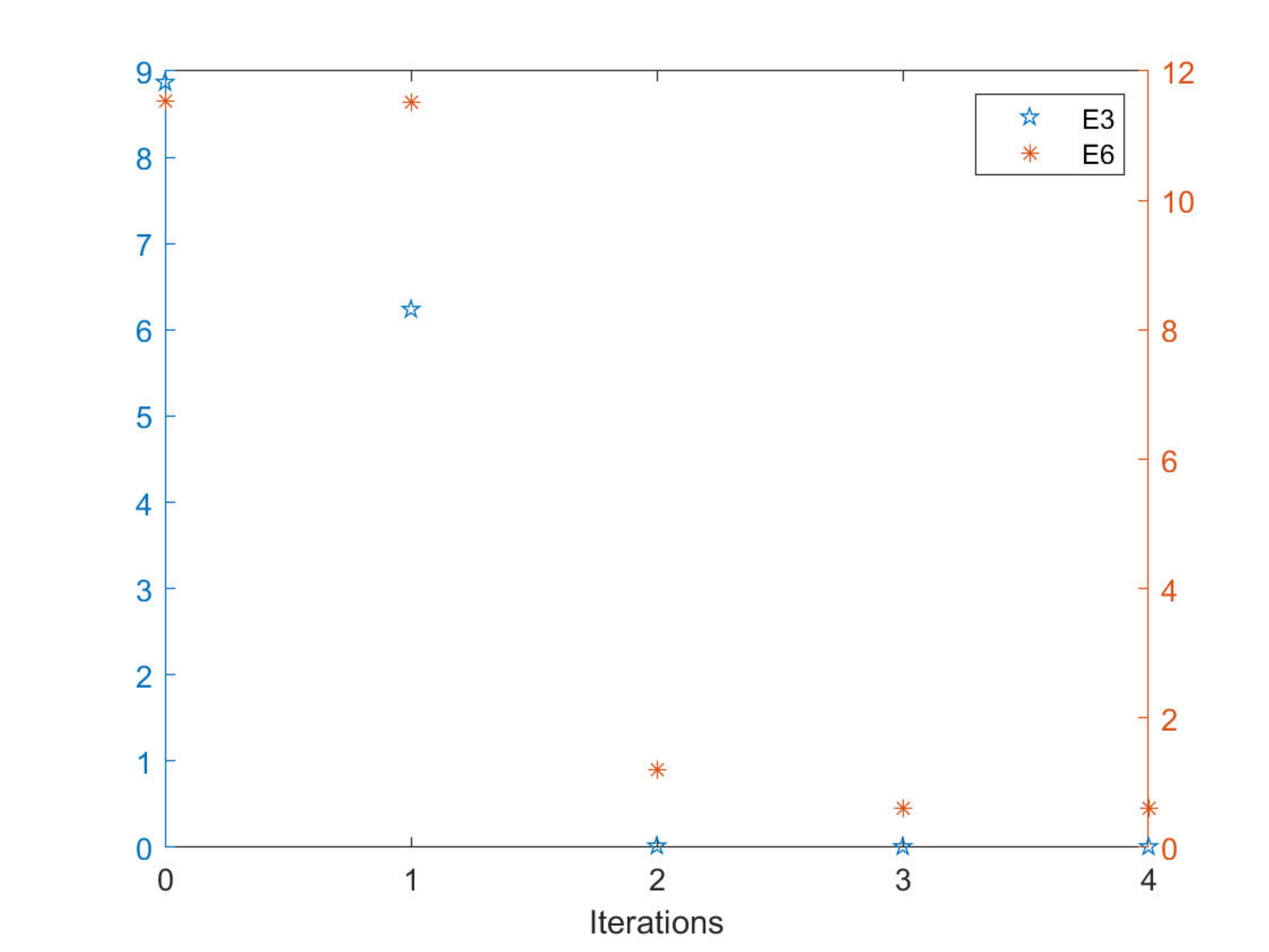}
\caption{\small{The evolution of the errors in the updating processes for the final selected $\kappa_1, \kappa_2$ and $\theta_T$. The left figure is for case I with $(\kappa_1, \kappa_2)=(0.375, 0.22)$ and $\theta_T=0.1$. The middle one is for case II with $(\kappa_1, \kappa_2)=(0, 0.01)$ and $\theta_T=0.2$. The right one is for case III with $(\kappa_1, \kappa_2)=(0.4, 0.3)$ and $\theta_T=0.075$. The update will be terminated at a setted 50 iterations if the terminal condition has not reached earlier. The last results are taken as the values in Figure \ref{fig2}.}}\label{fig1}
\end{center}
\end{figure}

  In this model reconstruction problem, taking $\kappa_1=0.025 k, \kappa_2=0.01 j$, for $k=0, 1, \cdots, 16$, $j=1, 2, \cdots, 30$. For a fixed $(\kappa_1, \kappa_2)$, starting from a given $\theta_T$, we constantly update $\vec{b}$ and solve $(\vec{\beta}, \varepsilon)$ from $\vec{b}$ until a stable estimator $(\vec{\beta}, \varepsilon)$ is obtained, and calculate the errors by the obtained estimators. The evolution of the errors in each updation is illustrated in Figure \ref{fig1} for case I, case II and case III, from left to right, respectively.

\begin{figure}[htbp]
\begin{center}
\includegraphics[width=5.3cm, height=6cm]{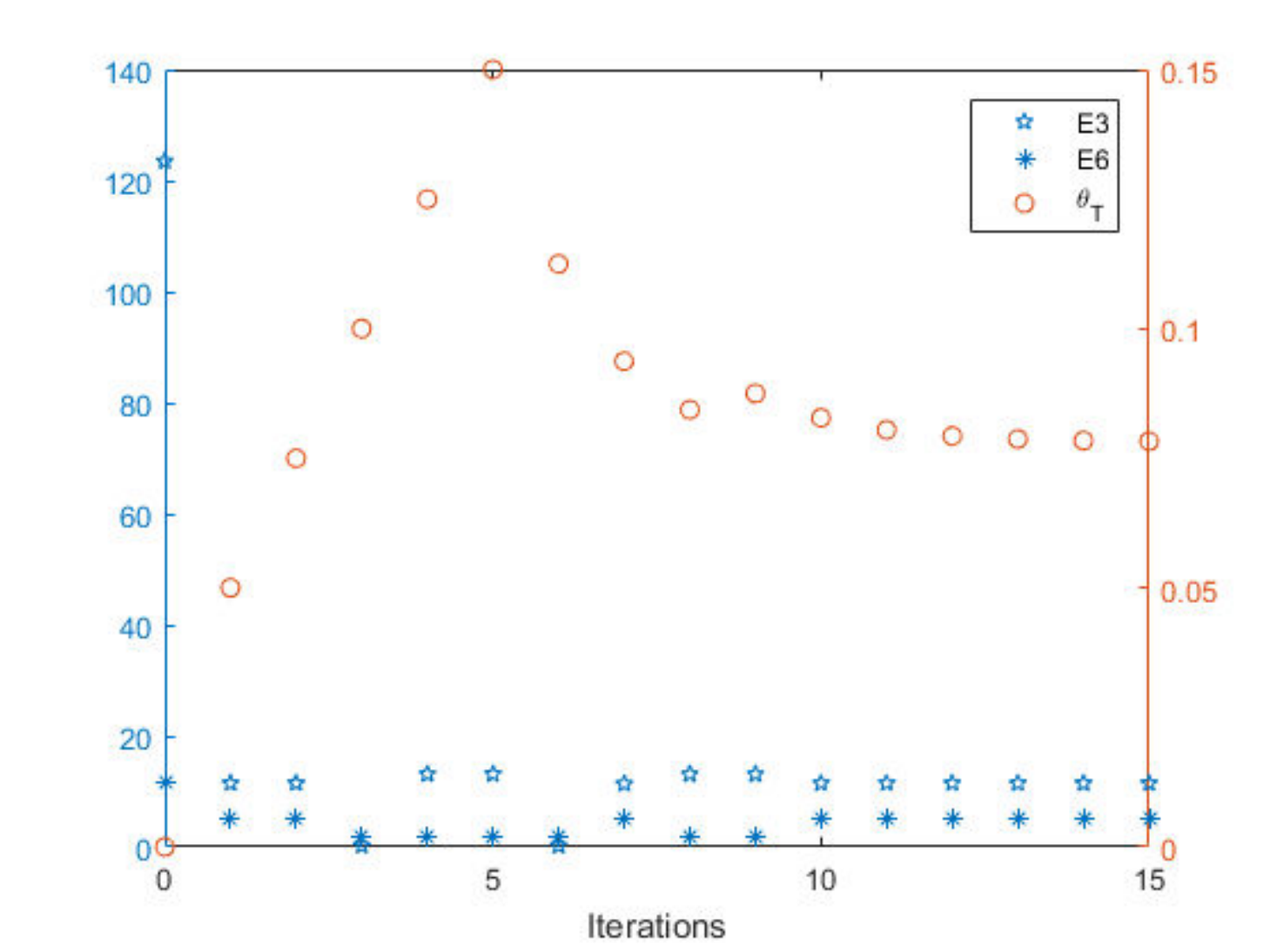}
\includegraphics[width=5.3cm, height=6cm]{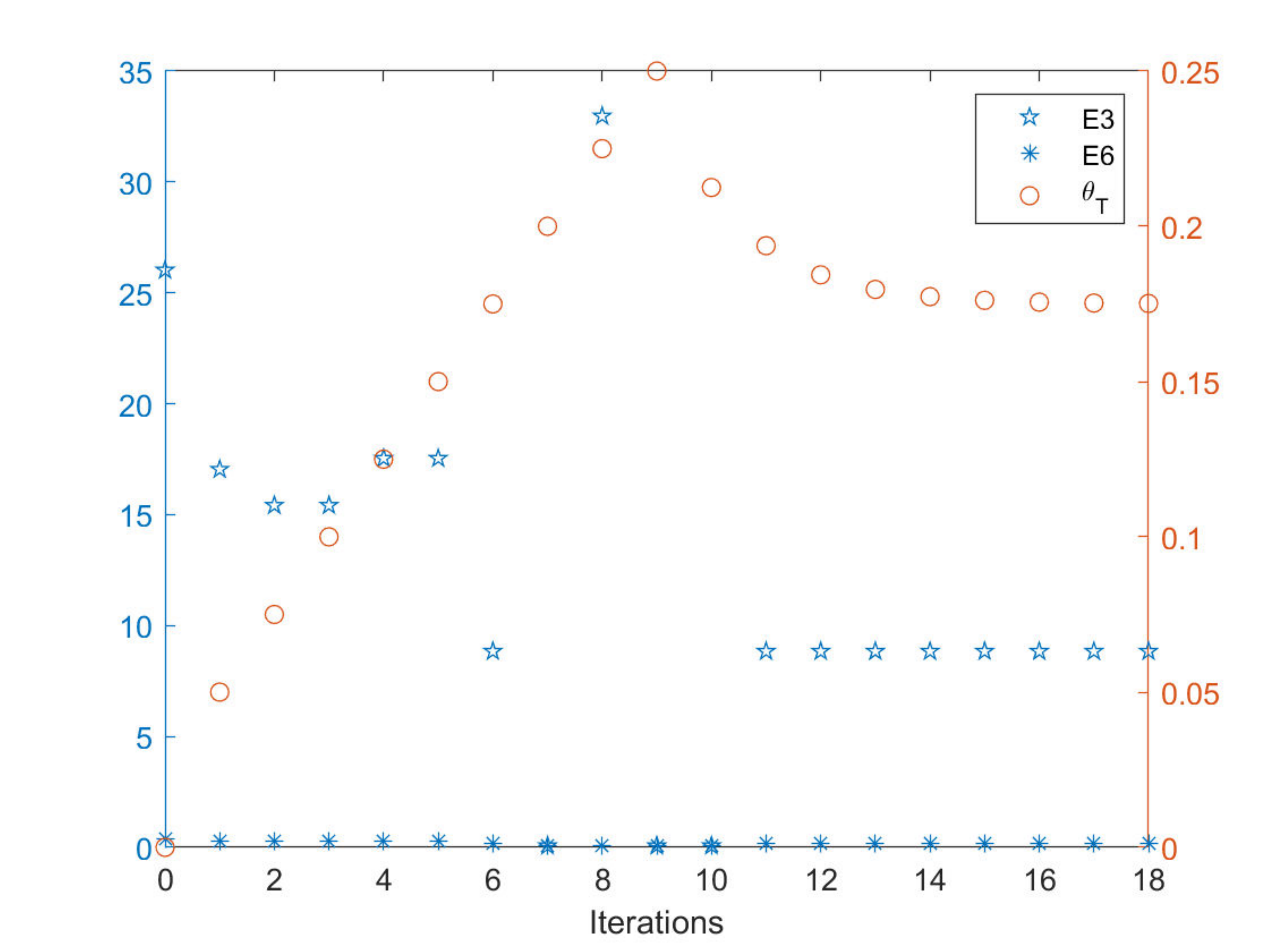}
\includegraphics[width=5.3cm, height=6cm]{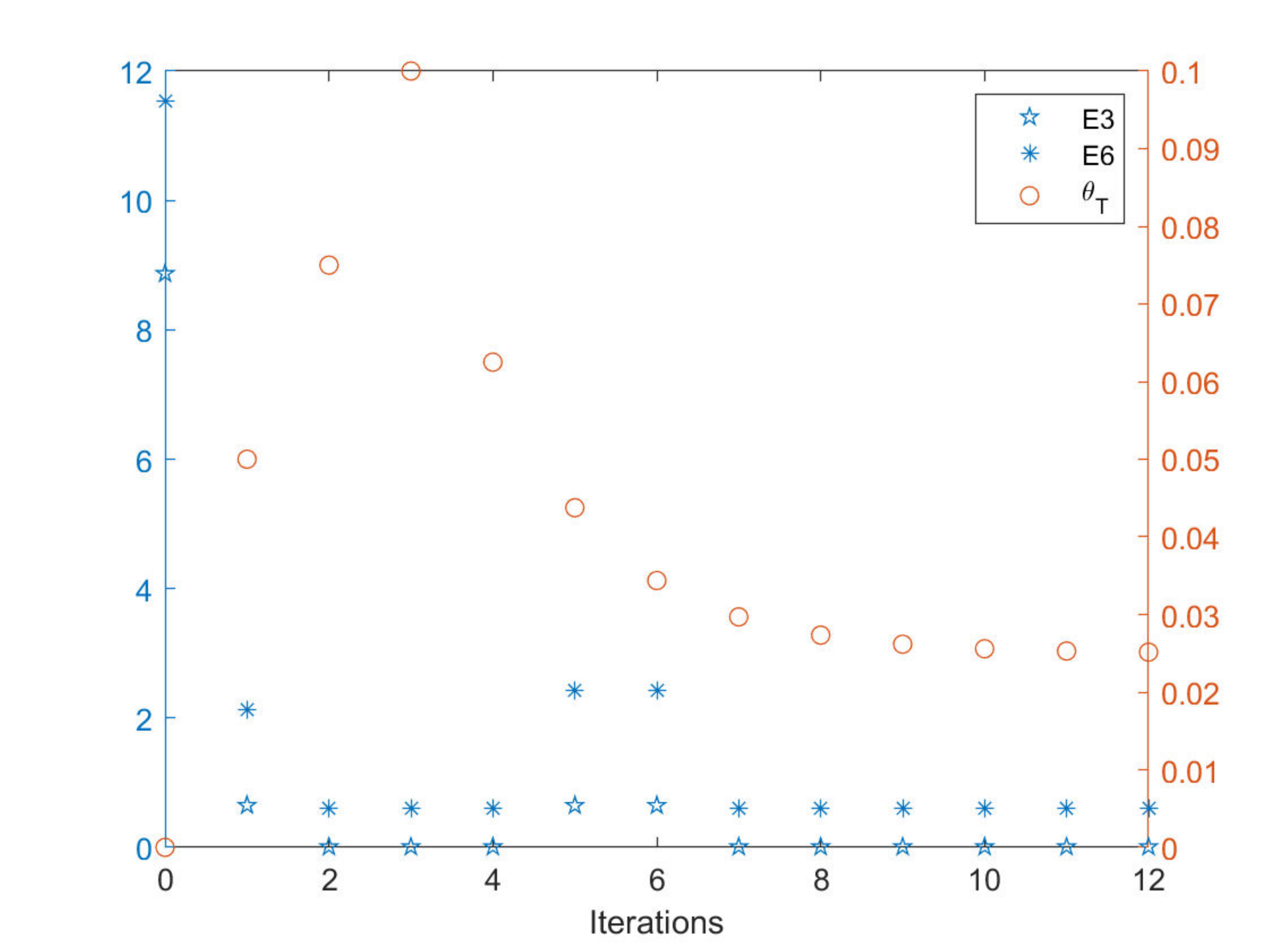}
\caption{\small{The errors after update for a threshold and the final selected $\kappa_1$ and $\kappa_2$. The three figures are for case I, II, and III, from left to right, with $(\kappa_1, \kappa_2)$ being $(0.375, 0.22)$, $(0, 0.01)$ and $(0.4, 0.3)$, respectively. The decrease of $E6$ enlarges the next threshold $\theta_T$, and in turn the increase of $E6$ pull back $\theta_T$.}}\label{fig2}
\end{center}
\end{figure}

  For a fixed $(\kappa_1, \kappa_2)$, we then compare the error $E6$ with the old best one to collect better estimator and regenerate $\theta_T$ by its renewal rule.  Repeat the above update process with the new $\theta_T$, until the new generated $\theta_T$ is close enough to the older one. With a threshold, the errors after update are shown in Figure \ref{fig2}.

%

\begin{figure}[htbp]
\begin{center}
\includegraphics[width=5.3cm, height=6cm]{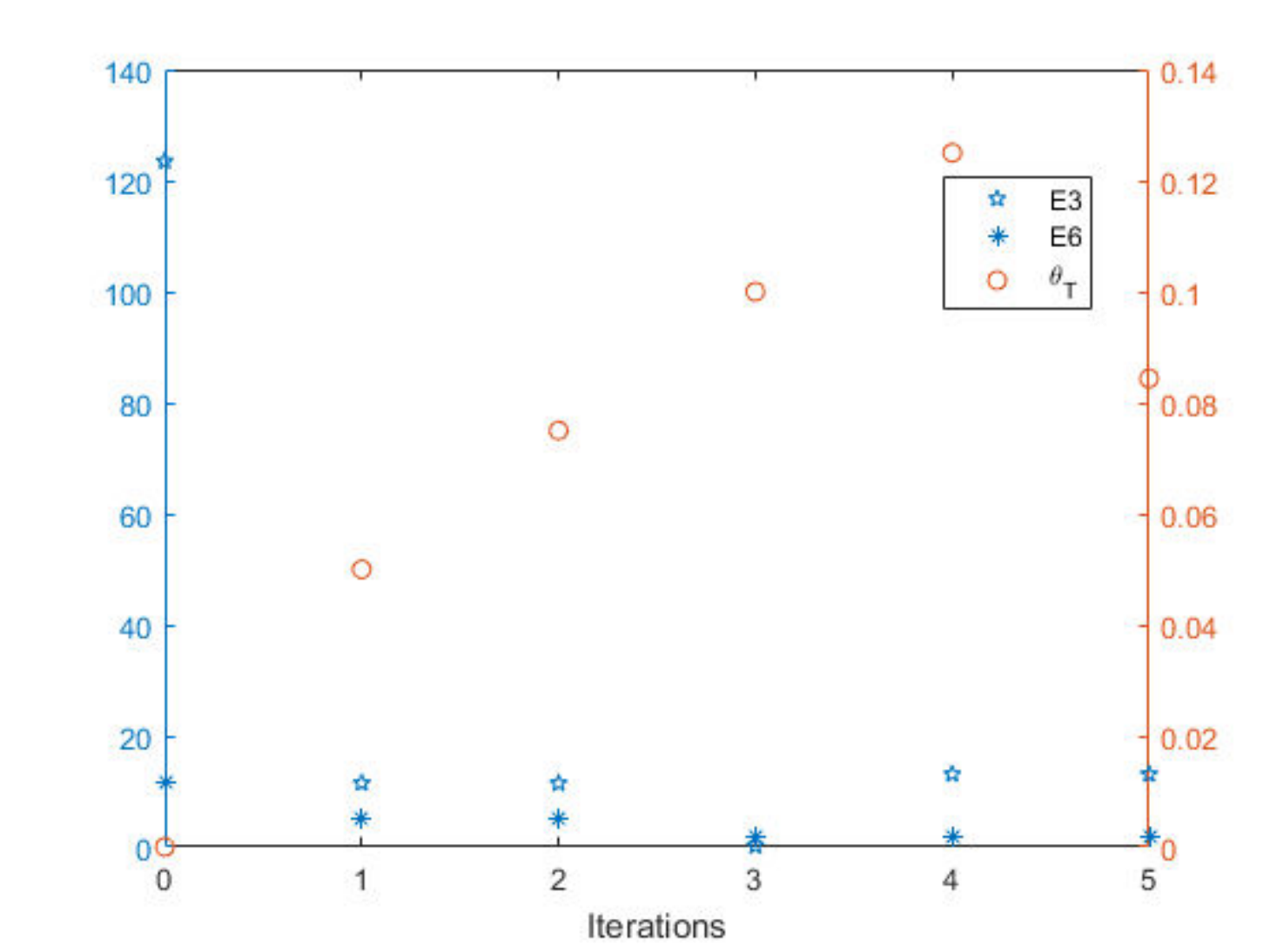}
\includegraphics[width=5.3cm, height=6cm]{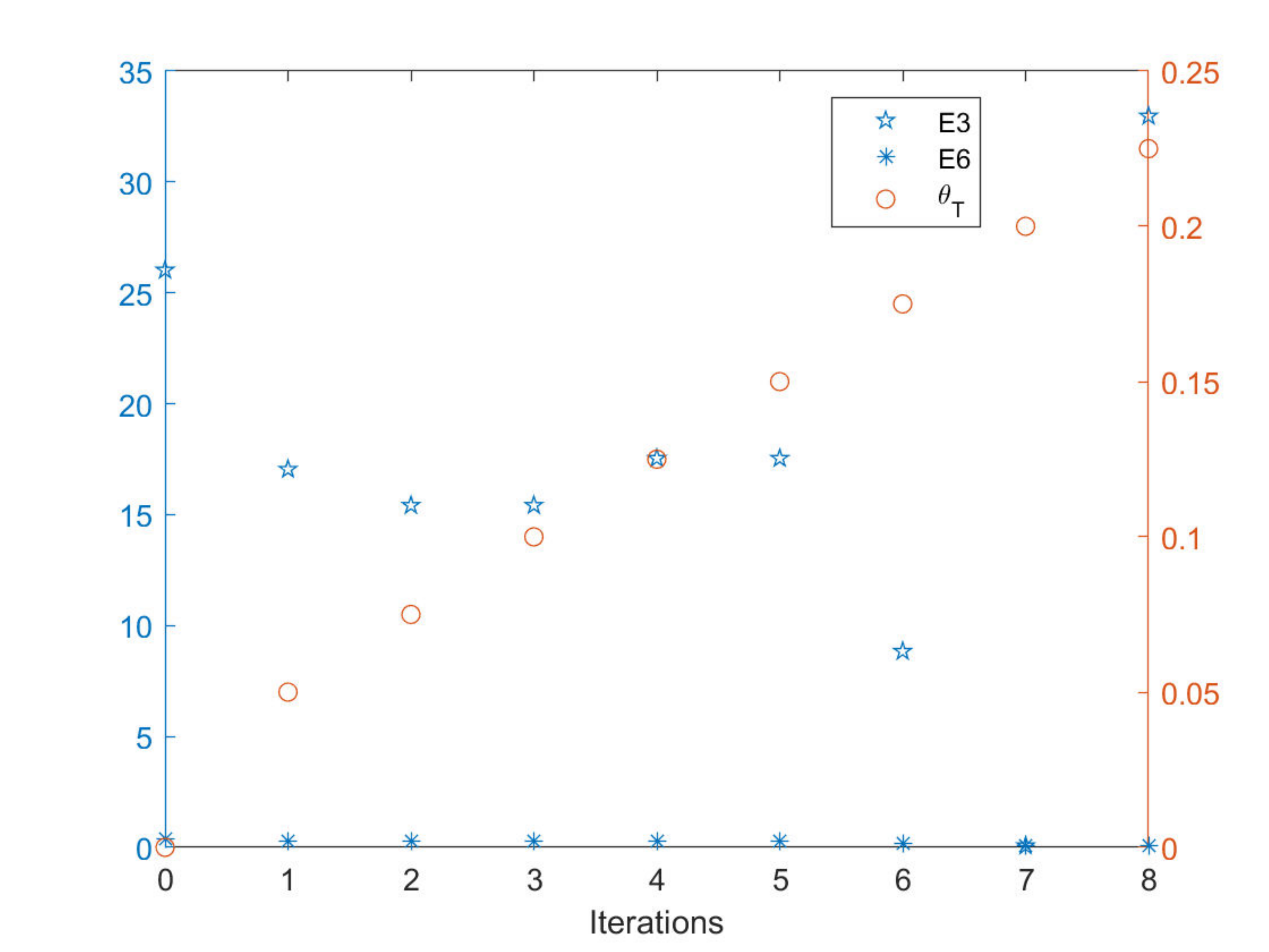}
\includegraphics[width=5.3cm, height=6cm]{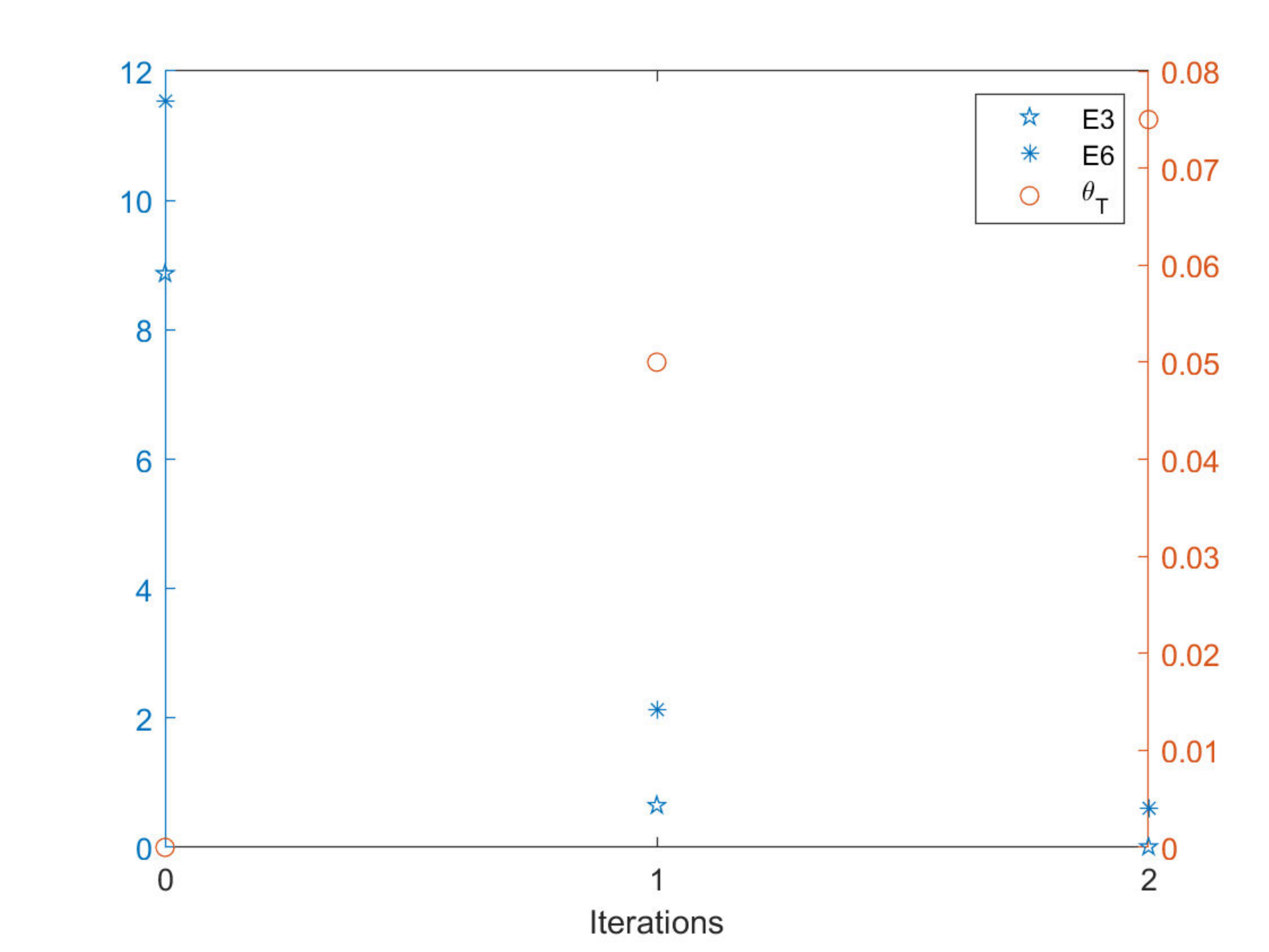}
\caption{\small{The selected best results. The three figures are for case I, II, and III, from left to right, with $(\kappa_1, \kappa_2)$ being $(0.375, 0.22)$, $(0, 0.01)$ and $(0.4, 0.3)$ respectively. The best results are selected from the updated results in Figure \ref{fig2} according to the principle of the error $E6$ turns smaller and smaller. }}\label{fig3}
\end{center}
\end{figure}

 We select the results with $E6$ in Figure \ref{fig2} smaller and smaller as the best results. The errors $E3, E6$ and $\theta_T$ of the best results are shown in Figure \ref{fig3}.  From these results, we then select the threshold with the smallest $E3$ as the threshold for the fixed $(\kappa_1, \kappa_2)$. For case I, II and III in Figure \ref{fig3}, $\theta_T$ equals 0.1, 0.2 and 0.075, respectively.

\begin{figure}[htbp]
\begin{center}
\includegraphics[width=15cm, height=11cm]{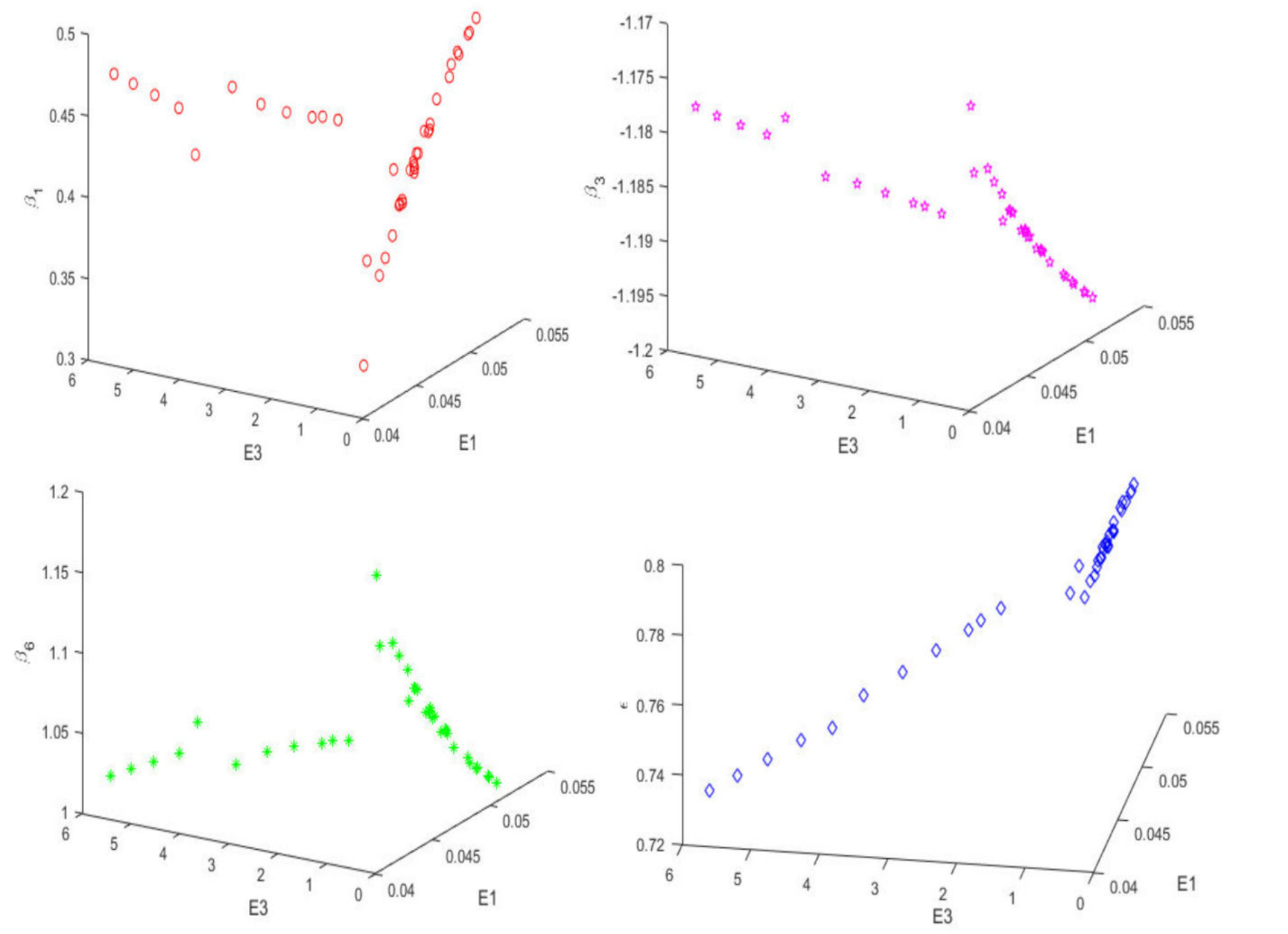}
\caption{\small{The errors $E1, E3$ v.s. $\beta_1$ (upper left), $\beta_3$ (upper right), $\beta_6$ (lower left) and $\varepsilon$ (lower right) for $(E4, E5)=(3, 8)$ listed in Table \ref{error11} and \ref{error12}.}}\label{fig4}
\end{center}
\end{figure}

After a threshold $\theta_T$ is selected, a pair of $(\kappa_1, \kappa_2)$ corresponds to a set of results. Sort by $E1$ from small to big, we listed the results with $E3<1$ in Table \ref{error11} and \ref{error12}. The errors with $E3\in(1, 2)$ are presented in the middle part of Table \ref{error12}. The errors with $E3>2$ and $(E4, E5)=(3, 8)$ are shown in the bottom of Table \ref{error12}.
We also listed the estimators with $(E4, E5)=(3, 8)$ in Table \ref{Betas}. These results are also sorted by $E1$ from small to big, with $E3<1$ in the upper part, and $E3>1$ in the half bottom part. This table together with the corresponding Figure \ref{fig4} illustrated that, when $E3<1$, $\beta_1$ and $\varepsilon$ both increase with $E1$ to a stable values, and $\beta_3, \beta_6$ each decreases with $E1$ to a stable values. When $E3$ exceeds 1, the values would go to the other side. Actually, we also observed that the estimators with $(E4, E5)$ take other values have no tendency to converge. So, we take the estimator with $(E1, E3)=(0.051521931, 0.231136499)$ and $(\kappa_1, \kappa_2)=(0.375, 0.22)$.

As listed in Table \ref{error2}, the minimum of $E3$ is 0.068423125 and then the value jumped to 2.414255893. So, we take the estimator for $(\kappa_1, \kappa_2)=(0, 0.01)$, which is listed in Table \ref{Beta}. In Table \ref{error3}, most of $E3<0.1$, and $(E4, E5)=(1, 2)$. In fact, all the estimators with $(E4, E5)=(1, 2)$ are close enough, and have a slight tendency with the increase of  $E1$.  We take the estimator corresponding to $(\kappa_1, \kappa_2)=(0.4, 0.3)$.

\newpage
\begin{table}[H]
\centering
\caption{Part of the selected results corresponding to $(\kappa_1, \kappa_2)$ for case I. }\label{error11}
\scalebox{0.95}{
\begin{tabular}{|c|c|c|c|c|c|c|c|}
  \hline
    & $\kappa_1$ & $\kappa_2$ & $E1$ & $E2$ & $E3$  &  $E4$ & $E5$  \\
  \hline
1&0.3&	0.12&	0.026939815&	0.056260831&	0.810529903&	4&	15 \\
2&0.35&	0.15&	0.027512134&	0.02881241&	0.403005153&	4&	15\\
3&0.275&0.07&	0.030155923&	0.089651519&	0.625205233&	6&	18\\
4&0.275&0.06&	0.031185227&	0.035798805&	0.3007625&	6&	18\\
5&0.35&	0.07&	0.032058741&	0.045525289&	0.33543262&	8&	30\\
6&0.375&0.06&	0.034403747& 0.063865974&	0.388980495&	8&	30\\
7&0.3&	0.05&	0.035396018&	0.009287062&	0.084640949&	8&	30\\
8&0.325&0.05&	0.035977481&	0.0304864&	0.18484621&	8&	30\\
9&0.25&	0.13&	0.040915848&	0.009159243&	0.186932482&	3&	8\\
10&0.225&0.18&	0.042568948&	0.012027928&	0.23459647&	3&	8\\
11&0.25&0.12&	0.042680852&	0.027285887&	0.533461953&	3&	8\\
12&0.3&	0.15&	0.042961344&	0.009841414&	0.19982028&	3&	8\\
13&0.35&0.21&	0.043641818&	0.009934937&	0.201417306&	3&	8\\
14&0.2&	0.15&	0.044437576&	0.012708587&	0.24357159&	3&	8\\
15&0.25&0.2&	0.044551965&	0.013315764&	0.256689297&	3&	8\\
16&0.375&0.2&	0.044584715&	0.010205395&	0.207388277&	3&	8\\
17&0.275&0.23&	0.044723258&	0.012662632&	0.244685253&	3&	8\\
18&0.25&0.21&0.045384028&0.031515335&0.587403886&	3&	8\\
19&0.1&	0.09&0.045432569&0.015458686&0.236955823&	3&	8\\
20&0.4&	0.21&0.045596116&0.009410949&0.192805735&	3&	8\\
21&0.325&0.16&0.045669816&0.009723&0.197954065	&3&	8\\
22&0.225&0.15&0.045744733&0.011669681&0.227695764&	3&	8\\
23&0.3&	0.25&0.045966062&0.015049886&0.288844346&	3&	8\\
24&0.225&0.17&0.046154114&0.014289434&0.272107785&	3&	8\\
25&0.275&0.22&0.046245681&0.013153949&0.25360326&	3&	8\\
26&0.125&0.12&0.046822365&0.016235436&0.255799042&	3&	8\\
27&0.375&0.19&0.046978399&0.009922707&0.202500177&	3&	8\\
28&0.4&	0.2&0.047077381&0.010054288&0.205377021&	3&	8\\
29&0.25&0.17&0.047244334&0.011755312&0.229606999&	3&	8\\
30&0.15&0.15&0.048045942&0.016967716&0.273986594&	3&	8\\
31&0.25&0.15&0.049011349&0.011463523&0.225690301&	3&	8\\
32&0.175&0.18&049449323&0.017466148&0.28680301&	3&	8\\
33&0.275&0.15&0.049906621&0.011569337&0.230190166&	3&	8\\
34&0.2&	0.2&0.050009585&0.016499411&0.282830571&	3&	8\\
35&0.325&0.19&0.050788088&0.011618899&0.232419407&	3&	8\\
36&0.325&0.18&0.050872958&0.01119349&0.22506516&	3&	8\\
37&0.375&0.22&0.051521931&0.011466825&0.231136499&	3&	8\\
38&0.125&0.11&0.052500087&0.031450179&0.68784801&	6&	16\\
39&0.375&0.21&0.061619699&0.0064827&0.131715685&	2&	6\\
40&0.4&	0.23&0.061631604&0.0060925&0.126164454&	2&	6\\
41&0.3&	0.11&0.091869213&0.056813276&0.722115448&	3&	6\\
42&0.375&0.14&0.091941327&0.057546552&0.708967114&	3&	6\\
43&0.35&0.13&0.135365079&0.073483969&0.657117275&	7&	20\\
44&0.25&0.1&0.193240972&4.44E-16&0.193240972&	2&	3\\
45&0.275&0.09&0.193501557&1.67E-16&0.193501557&	2&	3\\
\hline
\end{tabular}}
\end{table}

\begin{table}[H]
\centering
\caption{Part of the selected results corresponding to $(\kappa_1, \kappa_2)$ for case I. }\label{error12}
\scalebox{0.95}{
\begin{tabular}{|c|c|c|c|c|c|c|c|}
  \hline
    & $\kappa_1$ & $\kappa_2$ & $E1$ & $E2$ & $E3$  &  $E4$ & $E5$  \\
  \hline
46&0.3&	0.1&0.194227324&9.04E-16&0.194227324&	2&	3\\
47&0.325&0.1&0.194657739&3.33E-16&0.194657739&	2&	3\\
48&0.4&	0.1&0.196271553&1.67E-16&0.196271553&	2&	3\\
49&0.375&0.17&0.478744738&0.0022207&0.477349246&4&	6\\  \hline
50&0.35&0.14&0.024948326&0.129172171&1.88336774&5&19\\
51&0.275&0.08&0.025003415&0.190559162&1.533946863&7	&28\\
52&0.4&0.17&0.026388306&0.086913147&1.245602014&4&15\\
53&0.225&0.08&0.026447746&0.152348799&1.950616605&5&	17\\
54&0.375&0.16&0.026471506&0.081392076&1.170002948&4	&15\\
55&0.325&0.13&0.026573134&0.082975415&1.169601999&4	&15\\
56&0.25&0.07&0.027088397&0.12331372&1.734503499&7&	28\\
57&0.375&0.13&0.030524862&0.088563385&1.723396516&5&	17\\
58&0.35&0.09&0.031469022&0.102544481&1.495024936&7	&20\\
59&0.25&0.05&0.034614687&0.106404303&1.57749973&5	&18\\
60&0.175&0.01&0.036935561&0.054904564&1.915958055&5	&14\\
61&0.175&0.3&0.044739889&0.166985152&1.986306967&3	&8\\
62&0.15&0.24&0.044949533&0.143056378&1.70060192&3	&8\\
63&0.3&0.24&0.049352375&0.065626621&1.63626833&4	&11\\
64&0.275&0.21&0.049627393&0.067845621&1.731876756&4	&11\\
65&0.325&0.27&0.049965305&0.069480274&1.723258751&4	&11\\
66&0.325&0.26&0.050523389&0.071902797&1.794211573&4	&11\\
67&0.35&0.29&0.051215719&0.075963758&1.886780118&4	&11\\
68&0.35&0.28&0.051618921&0.0776454&1.938630281&4	&11\\
69&0.225&0.03&0.070470729&0.116532936&1.833947559&6	&16\\
70&0.35&0.06&0.09190076&0.062677789&1.054933993&3	&6\\
71&0.4&0.06&0.094183134&0.065509401&1.121334248&3	&6\\
72&0.35&0.11&0.155140184&0.104546119&1.91901288&7	&20\\
73&0.1&0.04&0.16118239&0.073680877&1.452335064&5	&16\\   \hline
74&0.05&0.09&0.040985874&0.364451385&3.880124457&3&8 \\
75&0.1&0.23&0.041070783&0.534474352&5.678033374&3&8 \\
76&0.1&0.22&0.041190925&0.493802891&5.286181444&3&8 \\
77&0.1&0.21&0.041345796&0.449789123&4.853110819&3&8 \\
78&0.1&0.2&0.041538487&0.402130029&4.374671742&3&8 \\
79&0.15&0.28&0.043135411&0.316835269&3.58159286	&3	&8 \\
80&0.15&0.27&0.0434484&0.264154135&3.028367111	&3	&8 \\
81&0.15&0.26&0.043873034&0.221421481&2.568209367&	3&	8 \\
82&0.15&0.25&0.044374308&0.181447426&2.12985506	&3	&8 \\
  .&.	&. &. &.	&.  & .  &. \\
  .&.	&.	&.  &.	&.  & .  &.  \\
  .&.	&.	&.  &.	&.  & .  &.  \\
  \hline
\end{tabular}}
\end{table}

\begin{table}[H]
\centering
\caption{The estimators before basis and $dB_t$ for case I with $(E4, E5)=(3, 8)$. }\label{Betas}
\scalebox{0.8}{
\begin{tabular}{|c|c|c|c|c|c|c|c|c|c|c|c|}
  \hline
   $(\kappa_1, \kappa_2)$& $1$ & $X$ & $X^2$ & $X^3$ & $X^4$ & $X^5$ & $\sin{X}$& $\cos{X}$ & $\sin{2X}$&$\cos{2X}$ & $dB_t$  \\
  \hline
  $(0.25,0.13)$ &0&0.327836541&0&-1.172722151&0&0&1.181206256&0&0&0&0.795855492\\
  $(0.225,0.18)$&0&0.376138461&0&-1.17957815&0&0&1.131442228&0&0&0&0.795463372\\
  $(0.25,0.12)$&0&0.382902527&0&-1.1803239&0&0&1.127508687&0&0&0&0.791033626\\
  $(0.3,0.15)$&0&0.385527258&0&-1.181034891&0&0&1.121778617&0&0&0&0.79591346\\
  $(0.35,0.21)$&0&0.396378089&0&-1.182581605&0&0&1.110499872&0&0&0&0.79596499\\
  $(0.2,0.15)$ &0&0.411486455&0&-1.184661332&0&0&1.095005472&0&0&0&0.795493667\\
  $(0.25,0.2)$ &0&0.411941703&0&-1.18471039&0&0&1.094619837&0&0&0&0.795315485\\
  $(0.375,0.2)$ &0&0.41256155&0&-1.184925465&0&0&1.093918174&0&0&0&0.795928136\\
  $(0.275,0.23)$&0&0.413602104&0&-1.184959797&0&0&1.092782024&0&0&0&0.79550349\\
  $(0.25,0.21)$&0&0.427591707&0&-1.186520317&0&0&1.080785909&0&0&0&0.790766547\\
 $(0.1,0.09)$&0&0.429209147&0&-1.18707282&0&0&1.075987734&0&0&0&0.79585546\\
$(0.4,0.21)$ &0&0.427307792&0&-1.187071538&0&0&1.078576652&0&0&0&0.796205009\\
$(0.325,0.16)$ &0&0.429485041&0&-1.187382523&0&0&1.076388386&0&0&0&0.79613473\\
$(0.225,0.15)$ &0&0.430746515&0&-1.187480369&0&0&1.075095588&0&0&0&0.7957778\\
$(0.3,0.25)$ &0&0.43165549&0&-1.187508313&0&0&1.074471091&0&0&0&0.794959944\\
$(0.225,0.17)$ &0&0.436294561&0&-1.188202355&0&0&1.06959379&0&0&0&0.795204042\\
$(0.275,0.22)$ &0&0.435975963&0&-1.188181238&0&0&1.069773658&0&0&0&0.795469545\\
$(0.125,0.12)$ &0&0.447286497&0&-1.189663626&0&0&1.057458227&0&0&0&0.79567531\\
$(0.375,0.19)$ &0&0.446381874&0&-1.189817485&0&0&1.058985891&0&0&0&0.796144045\\
$(0.4,0.2)$ &0&0.44750395&0&-1.189979512&0&0&1.057867511&0&0&0&0.796102671\\
$(0.25,0.17)$&0&0.450091126&0&-1.19026942&0&0&1.05513096&0&0&0&0.795840361\\
$(0.15,0.15)$ &0&0.461778709&0&-1.191737837&0&0&1.042632012&0&0&0&0.795485956\\
$(0.25,0.15)$ &0&0.471727846&0&-1.193411692&0&0&1.032837848&0&0&0&0.795971874\\
$(0.175,0.18)$ &0&0.47733876&0&-1.193971198&0&0&1.026658963&0&0&0&0.795378537\\
$(0.275,0.15)$ &0&0.481860578&0&-1.194885805&0&0&1.022493722&0&0&0&0.795929648\\
$(0.2,0.2)$  &0&0.482697535&0&-1.194782579&0&0&1.021245417&0&0&0&0.795416281\\
$(0.325,0.19)$ &0&0.490570278&0&-1.196139848&0&0&1.013507329&0&0&0&0.795941835\\
$(0.325,0.18)$ &0&0.491666267&0&-1.196313738&0&0&1.012351169&0&0&0&0.796042773\\
$(0.375,0.22)$ &0&0.497697853&0&-1.197175141&0&0&1.006158184&0&0&0&0.79598755\\
\hline
$(0.175,0.3)$&0&0.454174623&0&-1.186086693&0&0&1.050134233&0&0&0&0.775218978 \\
$(0.15,0.24)$&0&0.453119072&0&-1.186639601&0&0&1.051347485&0&0&0&0.778566674 \\
$(0.05,0.09)$&0&0.434770239&0&-1.177300571&0&0&1.065963916&0&0&0&0.753095468\\
$(0.1,0.23)$&0&0.473240359&0&-1.178016643&0&0&1.020894422&0&0&0&0.732799577\\
$(0.1,0.22)$&0&0.469057514&0&-1.178574206&0&0&1.026859337&0&0&0&0.737157871\\
$(0.1,0.21)$&0&0.464083461&0&-1.179124308&0&0&1.033573959&0&0&0&0.742033289\\
$(0.1,0.2)$&0&0.458266525&0&-1.179669721&0&0&1.041065253&0&0&0&0.747476787\\
$(0.15,0.28)$&0&0.469388964&0&-1.183798644&0&0&1.032344744&0&0&0&0.756506872\\
$(0.15,0.27)$&0&0.460914604&0&-1.184145608&0&0&1.042073815&0&0&0&0.762941264\\
$(0.15,0.26)$&0&0.456942839&0&-1.184845799&0&0&1.046777233&0&0&0&0.768328595\\
$(0.15,0.25)$&0&0.454510834&0&-1.185689335&0&0&1.049692376&0&0&0&0.773486361\\
  \hline
\end{tabular}}
\end{table}

\begin{table}[H]
\centering
\caption{Part of the selected results corresponding to $(\kappa_1, \kappa_2)$ for case II. }\label{error2}
\begin{tabular}{|c|c|c|c|c|c|c|c|}
  \hline
    & $\kappa_1$ & $\kappa_2$ & $E1$ & $E_2$ & $E_3$  &  $E4$ & $E5$  \\
  \hline
  1&0	&0.01&	0.068423125&	8.88E-16	&0.068423125	&2	&3  \\ \hline
  2&0.025&0.01&	0.003744129&	0.384764721	&2.414255893	&5	&20  \\
  3&0.05	&0.01&	0.006488225&	0.491778864	&6.5351293	    &7	&30  \\
  4&0.05	&0.02&	0.006719289&	0.530244365	&7.012510004	&6	&30  \\
  5&0.375&0.07&	0.092935253&	0.968517918	&7.066051677	&3	&12 \\
  6&0.4	&0.07&	0.095319797&	0.953445464	&7.163854265	&3	&12 \\
  7&0.075&0.01&	0.007170651&	0.489409673	&7.227866757	&7	&30 \\
  8&0.125&0.05&	0.008515818&	0.548073081	&7.5003069	    &6	&30 \\
  9&0.25	&0.08&	0.010016016&	0.553479956	&7.516204754	&6	&30 \\
  10&0.275&0.08&	0.010316231&	0.551248268	&7.532757737	&6	&30 \\
  11&0.2	&0.08&	0.009441848&	0.558228985	&7.536289451	&6	&30 \\
  .&.	&. &. &.	&.  & .  &. \\
  .&.	&.	&.  &.	&.  & .  &.  \\
  .&.	&.	&.  &.	&.  & .  &.  \\
  \hline
\end{tabular}
\end{table}

\begin{table}[H]
\centering
\caption{Part of the selected results corresponding to $(\kappa_1, \kappa_2)$ for case III. }\label{error3}
\begin{tabular}{|c|c|c|c|c|c|c|c|}
  \hline
    & $\kappa_1$ & $\kappa_2$ & $E1$ & $E_2$ & $E_3$  &  $E4$ & $E5$  \\
  \hline
1&0	&0.01	&0.00791007	&0&0.00791007&1	&2 \\
2&0	&0.02	&0.00791007	&0&0.00791007&1	&2 \\
3&0	&0.03	&0.007910072&0&0.007910072&1	&2 \\
4&0	&0.05	&0.007910082&5.55E-17&0.007910082&1	&2  \\
5&0	&0.06	&0.007910091&0&0.007910091&1	&2 \\
6&0	&0.14	&0.007910234&0	&0.007910234&1	&2  \\
7&0	&0.23	&0.007910546&1.24E-16&0.007910546&1	&2  \\
  .&.	&. &. &.	&.  & .  &. \\
  .&.	&.	&.  &.	&.  & .  &.  \\
  .&.	&.	&.  &.	&.  & .  &.  \\
440&0.4 &0.26	&0.00793866&1.24E-16	&0.00793866	&1	&2\\
441&0.4 &0.27	&0.00793884&0	&0.007938842&	1	&2\\
442&0.4 &0.28	&0.00793902&0	&0.007939025	&1	&2\\
443&0.4 &0.3	&0.00793904&0	&0.00793904	&1	&2\\
444&0.2 &0.03	&0.001043796&0.002244972	&0.064379145	&5&	8\\
445&0.025&0.05	&0.097883427&4.55E-15	&0.097883427	&3	&4\\
446&0.1 &0.02	&0.099026845&3.33E-16	&0.099026845	&3	&4\\
447&0.125&0.01	&0.099642975&1.46E-15	&0.099642975	&3	&4\\
448&0.125&0.05	&0.100840545&9.22E-16	&0.100840545	&3	&4\\
.&.	&. &. &.	&.  & .  &. \\
  .&.	&.	&.  &.	&.  & .  &.  \\
  .&.	&.	&.  &.	&.  & .  &.  \\
  \hline
\end{tabular}
\end{table}

\begin{table}[H]
\centering
\caption{The true values and their estimators for case I, II, and III in \eqref{Ex1}. }\label{Beta}
\begin{tabular}{|c|cc|cc|cc|}
  \hline
   &  True value I &  Estimator &  True value II &  Estimator &  True value III &  Estimator\\
  \hline
  $1$       &  0      &  0          &  0      &  0          &  0      &  0     \\
  $X$       &  0.5    &  0.497697853   &  0.5    &  0.495295   &  0      &  0      \\
  $X^2$     &  0      &  0          &  0      &  0          &  0      &  0    \\
  $X^3$     &  -1.2   &  -1.197175141&  -1.2   &  -1.19413   &  0      &  0    \\
  $X^4$     &  0      &  0          &  0      &  0          &  0      &  0     \\
  $X^5$     &  0      &  0          &  0      &  0          &  0      &  0     \\
    \hline
  $\sin{X}$  &  1     &  1.006158184&  0      &  0          &  0      &  0    \\
  $\cos{X}$  &  0     &  0          &  0      &  0          &  1      &  1.002284     \\
  $\sin{2X}$ &  0     &  0          &  0      &  0          &  0      &  0     \\
  $\cos{2X}$ &  0     &  0          &  0      &  0          &  0      &  0     \\   \hline
  $dB_t$   &  0.8     &  0.79598755 &  0.8    &  0.80161    &  0.8    &  0.79607     \\
   \hline
\end{tabular}
\end{table}
\begin{table}[H]
\centering
\caption{The test errors, weight values and the threshold corresponding to the estimators in Table \ref{Beta}. }\label{error}
\begin{tabular}{|c|c|c|c|c|c|c|}
  \hline
      & $E1$ & $E2$ & $E3$ & $\kappa_1$ & $\kappa_2$  &  $\theta_T$ \\
  \hline
  Case I &0.051521931 & 0.011466825  & 0.231136499 & 0.375  & 0.22   & 0.1 \\
  Case II & 0.068423 & 8.88$e^{-16}$  & 0.068423 & 0      &0.01    & 0.2 \\
  Case III & 0.0079    & 0           &  0.0079  & 0.4    & 0.3    & 0.075\\
  \hline
\end{tabular}
\end{table}

\section{Discussion} \label{discussion}

 Based on the   fact that the most probable transition trajectory of an SDE can be described by an Euler-Lagrange equation (a deterministic second order ODE), we have devised an algorithm to discover the drift term and the diffusion coefficient of the SDE, from a time-series data of the most probable transition trajectory.

 We have employed polynomials and trigonometric functions as the basis of the drift term. We considered a one-dimensional model example but with different drifts. For more complex drifts or high dimensional cases, our method   still works but becomes more complicated.

 For example,  for an SDE with a rational function in the drift,      a naive idea    is  to multiply    both sides  of the SDE  by the denominator of the rational function.  This converts the rational drift to a polynomial drift \cite{Mangan}. More precisely, for an SDE
 \begin{eqnarray*}
  dX=\Big(f(X)+\frac{g(X)}{G(X)}\Big) dt + \varepsilon dB_t,
 \end{eqnarray*}
the corresponding  Euler-Lagrange equation for the most probable transition trajectory $z$ is
  \begin{eqnarray*}
  \ddot z=\frac{\varepsilon^2}{2}\Big(f(z)+\frac{g(z)}{G(z)}\Big)'' + \Big(f(z)+\frac{g(z)}{G(z)}\Big)'\Big(f(z)+\frac{g(z)}{G(z)}\Big).
 \end{eqnarray*}
 Taking some derivatives, we get
 \begin{eqnarray*}
  G(z)^3\ddot z&=&\frac{\varepsilon^2}{2}\big(
  G^3f''+G^2g''-GG''g-2GG'g'+2G'G'g\big)\\ && + \big(G^3ff'+G^2f'g+G^2fg'+Ggg'-GG'fg-G'g^2\big).
 \end{eqnarray*}
 We write $G(z)^3=G_0 + G_1(z)$, where $G_0$ is a constant, and $G_1(z)$ is a polynomial without constant. The equation above is then
 \begin{eqnarray}
   G_0\ddot z&=&\Big\{\frac{\varepsilon^2}{2}\big(
  G^3f''+G^2g''-GG''g-2GG'g'+2G'G'g\big)\nonumber\\ && + \big(G^3ff'+G^2f'g+G^2fg'+Ggg'-GG'fg-G'g^2\big)\Big\}-\big[\ddot z G_1(z)\big].\label{rational}
 \end{eqnarray}
For another example, consider    a  two dimensional SDE system
 \begin{eqnarray*}
 dX_1 &=& f_1(X_1, X_2) dt + \varepsilon_1 dB_t^1, \\
 dX_2 &=& f_2(X_1, X_2) dt + \varepsilon_2 dB_t^2,  \\
 X(0) &=& \vec{a}_0\in \mathbb{R}^2,  \quad X(1) = \vec{b} \in \mathbb{R}^2.
 \end{eqnarray*}
 In searching for the most probable transition trajectory $z=(z_1, z_2)^T$, the corresponding OM function is
 \begin{eqnarray*}
 OM(z, \dot z)=\big(\frac{\dot z_1-f_1(z_1, z_2)}{\varepsilon_1}\big)^2+
\big(\frac{\dot z_2-f_2(z_1, z_2)}{\varepsilon_2}\big)^2+\frac{\partial f_1}{\partial z_1}+\frac{\partial f_2}{\partial z_2}.
 \end{eqnarray*}
 Then by the variation principle, $z(t)$ satisfies the Euler-Lagrange equation
 \begin{eqnarray}
 \ddot{z}_1 &=& \Big\{f_1\frac{\partial f_1}{\partial z_1} + \frac{\varepsilon_1^2}{\varepsilon_2^2}f_2\frac{\partial f_2}{\partial z_1}+ \frac{\varepsilon_1^2}{2}\big(\frac{\partial^2 f_1}{\partial z_1^2}+\frac{\partial^2 f_2}{\partial z_2\partial z_1}\big)\Big\} + \Big[\dot z_2\big(\frac{\partial f_1}{\partial z_2}-\frac{\varepsilon_1^2}{\varepsilon_2^2}\frac{\partial f_2}{\partial z_1}\big)\Big], \nonumber \\
\ddot{z}_2 &=& \Big\{f_2\frac{\partial f_2}{\partial z_2} + \frac{\varepsilon_2^2}{\varepsilon_1^2}f_1\frac{\partial f_1}{\partial z_2}+ \frac{\varepsilon_2^2}{2}\big(\frac{\partial^2 f_2}{\partial z_2^2}+\frac{\partial^2 f_1}{\partial z_1\partial z_2}\big)\Big\} + \Big[\dot z_1\big(\frac{\partial f_2}{\partial z_1}-\frac{\varepsilon_2^2}{\varepsilon_1^2}\frac{\partial f_1}{\partial z_2}\big) \Big]. \label{2dim}
 \end{eqnarray}
It is similar for both \eqref{rational} and \eqref{2dim}, when  we express each term in the drift in a polynomial basis. The part in the curly braces on the right hand side  of them are the expressions of polynomials, and the part in the square brackets are  the expressions of polynomials multiplying $\ddot z$, $\dot z_1$ or $\dot z_2$.  Thus, the corresponding expressions similar to $\vec{b}$ in \eqref{b_beta1}-\eqref{b_beta3} are    much more  complicated here,  and this increases the computational cost. 


\section*{Acknowledgements}
This work was partially supported by NSFC grants 11531006, 11601491 and 11771449.

\end{document}